\documentclass[reqno]{amsart}

\usepackage{geometry}
\usepackage{hyperref}
\usepackage[english]{babel}
\usepackage{latexsym, mathrsfs}
\usepackage{amsfonts, amssymb}
\usepackage{amsaddr}

\newtheorem{prop}{Proposition}[section]
\newtheorem{lem}[prop]{Lemma}
\newtheorem{cor}[prop]{Corollary}
\newtheorem{thm}[prop]{Theorem}
\newtheorem{con}[prop]{Conjecture}

\theoremstyle{definition}
\newtheorem{oss}[prop]{Remark}
\newtheorem{defn}[prop]{Definition}
\newtheorem{ex}[prop]{Example}

\usepackage[utf8]{inputenc}

\newcommand{\Z}{\mathbb{Z}}

\newcommand{\E}{\mathcal{E}}

\newcommand{\equad}{\quad \text{ and } \quad}
\newcommand{\lcm}{\mathrm{lcm}}

\def\MR{\mathrm{MR}}
\def\MRS{\mathrm{MRS}}

\numberwithin{equation}{section}

\usepackage[table,dvipsnames]{xcolor}

\hypersetup{ 
 colorlinks =true,
 linkcolor=OliveGreen,
 citecolor=Plum,
 urlcolor=olive,
 }

\begin{document}

\title[On a conjecture by
Sylwia Cichacz and Tomasz Hinc,
and a related problem]{On a conjecture by
Sylwia Cichacz and Tomasz Hinc,\\
and a related problem}

\author[Fiorenza Morini]{Fiorenza Morini}
\address{Dipartimento di Scienze Matematiche, Fisiche e Informatiche, Universit\`a di Parma,\\
Parco Area delle Scienze 53/A, 43124 Parma, Italy}
\email{fiorenza.morini@unipr.it}

\author[Marco Antonio Pellegrini]{Marco Antonio Pellegrini}
\address{Dipartimento di Matematica e Fisica, Universit\`a Cattolica del Sacro Cuore,\\
Via della Garzetta 48, 25133 Brescia, Italy}
\email{marcoantonio.pellegrini@unicatt.it}

\author[Stefania Sora]{Stefania Sora}
\address{Universit\`a Cattolica del Sacro Cuore,\\
Via della Garzetta 48, 25133 Brescia, Italy}
\email{stefania.sora01@icatt.it}

\begin{abstract}
A $\Gamma$-magic rectangle set
    $\MRS_\Gamma (a, b; c)$ is a collection of $c$ arrays  of size $a\times b$ whose entries are the elements of an abelian group $\Gamma$ 
    of order $abc$, each one appearing once and in a unique array in such a way that the sum of the elements of each row is equal to a constant $\omega \in \Gamma$ and the sum of the elements of each column is equal to a constant $\delta \in \Gamma$. 
    
In this paper we provide 
new evidences
for the validity of a conjecture proposed by Sylwia Cichacz and Tomasz Hinc on the existence of an $\MRS_\Gamma(a,b;c)$.
We also generalize 
this problem, describing constructions of $\Gamma$-magic rectangle sets, whose elements are partially filled arrays.
\end{abstract}

\maketitle

\section{Introduction}

Magic squares are some of the oldest combinatorial objects that have always fascinated mathematicians.
In this paper we consider one of the many generalizations that have been recently proposed.
First of all, we recall what a magic rectangle is.

\begin{defn}\label{MR}
    A \textit{magic rectangle} $\MR(a, b)$ is an  $a\times b$ array whose entries are the integers $1,2,\ldots,ab$, each one appearing exactly once in such a way that the elements of each row add up to a constant $\omega$ and the elements of each column add up to a constant $\delta$.
\end{defn}

\begin{thm}\cite{Ha1,Ha2}
    \label{th:esistenzarettmagico}
   A magic rectangle  $\MR(a, b)$ exists if and only if  $a,b>1$, $ab>4$ and $a\equiv b \pmod 2$. 
\end{thm}

Dalibor Froncek introduced in \cite{F13} the notion of magic rectangle set, also providing necessary and sufficient conditions for  its existence.

\begin{defn}\label{defMRS}
A \textit{magic rectangle set} $\MRS(a, b; c)$ is
a collection of  $c$  arrays of size $a\times b$, with entries from $\{1,2,\ldots,abc\}$, each appearing once and in a unique array,
with all row sums in each array equal to a constant $\omega$ and all column sums in each array equal to a constant $\delta$.
\end{defn}

\begin{thm}\cite{F13,F17}\label{th:Fr}
    Let $a,b>1$. An $\MRS(a, b; c)$ exists if and only if one of the following two conditions holds:
    \begin{itemize}
        \item[(1)] $a$ and $b$ are both even;
        \item[(2)] the product $abc$  is odd.
    \end{itemize}  
\end{thm}

Sylwia Cichacz and Tomasz Hinc introduced in \cite{CH21} the concept of a $\Gamma$-magic rectangle set.

\begin{defn}\label{GMRS}
    A \textit{$\Gamma$-magic rectangle set} 
    $\MRS_\Gamma (a, b; c)$ is a collection of $c$ arrays of size $a\times b$ whose entries are the elements of an abelian group $(\Gamma,+)$ 
    of order $abc$, each one appearing once and in a unique array in such a way that the sum of the elements of each row is equal to a constant $\omega \in \Gamma$ and the sum of the elements of each column is equal to a constant $\delta \in \Gamma$. 
\end{defn}

Considering the existence of these objects, Cichacz and Hinc proved the following result, where $\Upsilon$ denotes the set of all finite abelian groups of odd order or containing more than one involution.

\begin{thm}\cite[Theorem 3.9]{CH21b}\label{th:esistenza}
    Let $a,b>1$ be such that
    $\{a,b\}\neq\{2\ell+1,2^\alpha\}$ where $\alpha,\ell>0$. An $\MRS_\Gamma(a,b;c)$ exists if and only if either $a$ and $b$ are both even or $\Gamma\in\Upsilon$.
\end{thm}

They also proposed the following.

\begin{con}\cite{CH21b}\label{ConCH}
Let $a,b>1$. An $\MRS_\Gamma(a,b;c)$ exists if and only if either $a$ and $b$ are both even, or  $\Gamma\in\Upsilon$ and $\{a,b\}\neq\{2\ell+1,2\}$.
    \end{con}

Note that, to prove the previous conjecture, one has to consider only the case
$\{a,b\}=\{ 2\ell+1, 2^\alpha\}$,
where $\alpha, \ell>0$. 
In fact, in \cite{CH21} the same authors proved that a necessary and sufficient condition
for the existence of an $\MRS_\Gamma(2 \ell + 1, 4; 4 h + 2)$ is that the abelian group $\Gamma$ has more than one involution.

In this paper, we provide further evidences for the validity of Conjecture \ref{ConCH},
proving the following.

    \begin{thm}\label{main}
        An $\MRS_\Gamma(2\ell+1, 8; 4h+2)$ exists if and only if the abelian group $\Gamma$, of order $16(2\ell+1)(2h+1)$,  has more than one involution.
    \end{thm}

The last section of the paper is devoted to  constructions of $\Gamma$-magic rectangle sets whose elements are partially filled arrays, that is, arrays
where some cells are allowed to be empty. 

\section{Notation and preliminary results}

Given two integers $a,b$ such that  $a\leqslant b$, we denote by $[a,b]$ the set
$\{a,a+1,\ldots,b \}\subset \Z$.
If $a>b$, we set $[a,b]=\varnothing$.
Given an array (or, more in general, a partially filled array) $A$, we write $\E(A)$
for the list of its entries (considering their multiplicities). Similarly, if $T=\{A_1,A_2,\ldots,A_n\}$ is a
set of partially filled arrays, we write
$\E(T)$ for
$\E(A_1)\cup \ldots \cup\E(A_n)$.
If the sum of the entries in each row and in each column of $A$ is zero, then we say that $A$ is a \textit{zero-sum} array. Similarly for the set $T$, if this holds for every element of $T$. 
\smallskip

The elements of a finite cyclic group $(\Z_n,+)$ of order $n$ will be denoted by $[x]_n$. 
The fundamental theorem of finite abelian groups
says that any  nontrivial abelian group $\Gamma$ of order $n$ can be written as a direct 
sum of cyclic subgroups:
$$\Gamma\cong \Z_{d_1}\oplus \Z_{d_2}\oplus \ldots\oplus \Z_{d_s},$$
where $d_i\geqslant 2$ for all $i\in [1,s]$ and 
$d_i\mid d_{i+1}$ for all $i\in [1,s-1]$.
Note that $d_s$ is the exponent of $\Gamma$.
Furthermore, we can decompose $\Gamma$ as a direct sum of cyclic subgroups of prime power order:
$$\Gamma \cong \mathbb{Z}_{p_1^{\alpha_1}} \oplus \mathbb{Z}_{p_2^{\alpha_2}} \oplus \ldots \oplus \mathbb{Z}_{p_k^{\alpha_k}},$$
where $p_1,p_2,\ldots,p_k$ are (not necessarily distinct) primes and  $\alpha_i\geqslant 1$ 
for all $i\in [1,k]$. Clearly, one has
$$n=\prod_{i=1}^s d_i=\prod_{i=1}^k p_i^{\alpha_i}. $$

We recall (see, for instance,  \cite{S18})
that,  given a finite abelian group $\Gamma$,
the sum of its elements is $0_\Gamma$, unless
$\Gamma$ contains a unique involution $\iota$:
in this case, the sum of elements of $\Gamma$ is
$\iota$.
Furthermore, if $\Gamma$ has even order
and we decompose $\Gamma$ as a direct sum of cyclic subgroups, then
$\Gamma$ has exactly  $2^r-1$ involutions, where 
$r$ is the number of cyclic components 
of even order appearing in the decomposition.
Recall that  $\Upsilon$ denotes the set of all finite abelian groups having zero or at least two involutions.
We denote the order of an element  $g\in \Gamma$ by~$|g|$.
\smallskip

We now describe a possible connection between magic rectangle sets and distance magic labelings. Recall that a (simple and nonoriented) graph  $G$ is a pair $(V,E)$, 
where the elements of $V$ are called  \textit{vertices} of the graph and the elements of $E$ are unordered pairs of distinct elements of $V$, called  \textit{edges}.
The cardinality of  $V$ is called the  \emph{order} of the graph $G$.
Given a vertex $v$ of a graph $G$, we say that a vertex $w\in V$ is \textit{adjacent} to $v$ if the pair $\{v,w\}$ is an edge of $G$.
So, we can define the \textit{open neighborhood} of $v$ as the set $N(v)$ of the vertices of $V$ adjacent to $v$.
The \textit{degree} of a vertex $v\in V$ is the cardinality of $N(v)$.

A graph is said to be \textit{regular} if all its vertices have the same degree. If such degree is equal to a constant $k$, we say that the graph is $k$-regular. 
We recall that a $k$-regular graph of order $n$ has exactly $\frac{nk}{2}$ edges: hence, if  $k$ is odd, then $n$ must be even.
A graph $G$ is said to be \textit{$r$-partite} (where $r\geqslant 2$ is an integer) if the set  $V$ of the vertices admits a partition into $r$ classes in such a way that each edge of $G$ connects vertices of  distinct classes. 
An $r$-partite graph is complete if any two vertices belonging to distinct classes are connected by an edge.

\begin{defn}
    A \textit{distance magic labeling} of a graph $G=(V,E)$ of order $n$ is a bijection
    $\eta: V \to \{1,2,\ldots,n\}$ such that 
    the weight
    $$\omega(v)=\displaystyle \sum_{u\in N(v)} \eta(u)$$
    is equal to a constant  $\mu\in \Z$ for all $v\in V$.
    This integer $\mu$ is called the \textit{magic constant} of the labeling~$\eta$. 
\end{defn}

A graph that admits a distance magic labeling is called \textit{distance magic}. 
The interest about distance magic labeling is justified by the construction of magic rectangles, since we can build a distance magic complete 
$b$-partite graph, with each part of cardinality $a$, labeling the vertices of each part with the entries of the columns of a magic rectangle $\MR(a,b)$. 
The magic constant of  this graph is  $\mu=(b-1) \delta$, where $\delta$ is the sum of the elements of each column of the magic rectangle.

\begin{ex}\label{magicgr}
Starting from the following $\MR(3,5)$
    $$\begin{array}{|c|c|c|c|c|}\hline
   15 & 2 & 14 & 4 & 5 \\ \hline
   8  & 10 & 7 & 9 & 6\\ \hline
   1 &  12 & 3  & 11 & 13\\ \hline
    \end{array},$$
   one can obtain a  distance magic complete 
    $5$-partite graph, where each part consists of $3$ vertices. In this case, the magic constant is $\mu=96$.
    Clearly, considering the transpose of the previous array, one can obtain 
    a distance magic complete $3$-partite graph, where each part consists of $5$ vertices and the magic constant is equal to $80$.
\end{ex}

Froncek introduced in \cite{F13b} the notion of distance magic graph on an abelian group.

\begin{defn}
    A \textit{$\Gamma$-distance magic labeling} of a graph $G=(V,E)$, with $|V|=n$, is a bijection $\eta: V \to \Gamma$, where $(\Gamma,+)$ is an abelian group of order $n$, such that the weight
    $$\omega(x)=\displaystyle\sum_{y\in N(x)} \eta(y)$$  
    of each vertex $x\in V$ is equal to the same element $\mu \in \Gamma$, called the
    \textit{magic constant}.
\end{defn} 

A graph  $G=(V,E)$ is called \textit{group distance magic} if there exists a $\Gamma$-distance magic labeling for some abelian group $\Gamma$ of order $|V|$.  Similarly to what we did in Example
\ref{magicgr}, one can construct a 
$\Gamma$-distance magic complete $bc$-partite
graph, with each part of cardinality $a$, 
labeling the vertices of each part with the entries of the columns of every array in an
$\MRS_\Gamma(a,b;c)$.

We recall the following result.

\begin{thm}\cite[Theorem 19]{CF16}
    \label{th:graforegolare}
    Let $G$ be an $r$-regular graph on  $n$ vertices, where $r$ is odd. Then, there is no abelian group $\Gamma$ of order $n$ having a unique involution such that $G$ is
    $\Gamma$-distance magic.
\end{thm}

\section{Magic rectangle sets and generalizations}

As recalled in the Introduction, Froncek 
considered the existence of a magic rectangle set. On the other hand, 
Abdollah Khodkar and David Leach worked in \cite{KL} 
with magic rectangles where some cells are allowed to be empty.
So, as in \cite{MP4}, we are interested in the following objects.

\begin{defn}\label{MRS}
Let $m,n,s,k,c$ be five positive integers.
A magic rectangle set $$\MRS(m,n;s,k;c)$$ is a set
of $c$ partially filled $m\times n$ arrays
with entries in $[1,nkc]$ such that
    \begin{itemize}
        \item[{\rm (a)}] every entry
        appears  exactly once and in a unique array;
        \item[{\rm (b)}] for every array, each row contains exactly $s$ filled cells and each column contains exactly $k$ filled cells;
        \item[{\rm (c)}] for every array, the sum of the entries of each row is equal to  a constant  $\omega$ and the sum of the entries of each column is equal to a constant $\delta$.
    \end{itemize}
\end{defn}

If $nkc>1$, trivial necessary conditions for the existence of an $\MRS(m,n;s,k;c)$ are 
\begin{equation}\label{nec}
2\leqslant  s \leqslant n,\quad 2 \leqslant  k \leqslant m  \equad 
ms = nk. 
\end{equation}
Furthermore, note that we necessarily have 
$$\omega=\frac{s(nkc+1)}{2}\equad 
\delta=\frac{k(nkc+1)}{2}.$$

In \cite{MP4} it was proved that, when conditions \eqref{nec} are satisfied,
an $\MRS(m, n; s, k; c)$ exists if and only if either $nkc$ is odd, or $s$ and $k$ are both even and $sk > 4$.
This result relies on Theorem \ref{th:Fr}
that follows from a series of lemmas.
In particular, in \cite{F13} the following result was proved, but its proof has some flaws. For this reason, we give here 
a different argument.

\begin{lem}\cite[Lemma 4.3]{F13}\label{lem:MRS(2,b;c)}
    If $b\geqslant 4$ is even, then there exists an  $\MRS(2,b;c)$ for all $c\geqslant 1$.
\end{lem}

\begin{proof}
For every $i$, let
$$Q_i=\begin{array}{|c|c|c|c|} \hline
4i+1 & 2bc-4i-1 & 2bc-4i-2 & 4i+4 \\ \hline
2bc-4i & 4i+2 & 4i+3 & 2bc-4i-3 \\ \hline
\end{array}.$$
The sum of the elements in each row of $Q_i$ is $4bc+2$, while the sum of the elements in each column is $\delta=2bc+1$.
Furthermore, $\E(Q_i)=[4i+1,4i+4]\cup
[2bc-4i-3, 2bc-4i]$.

First, suppose that $b\equiv 0\pmod 4$.
In this case, we take
$$\mathcal{Q}_0=\left\{Q_i : i\in \left[0, \frac{bc}{4}-1\right] \right\}.$$
Since
$\E(\mathcal{Q}_0)=[1, bc]\cup 
[bc+1, 2bc]=[1,2bc]$, to obtain the elements of an $\MRS(2,b;c)$ it suffices to juxtapose $\frac{b}{4}$ distinct elements of $\mathcal{Q}_0$ constructing a $b\times 2$ rectangle and repeat this operation  $c$ times.
The sum of the elements of each row of these rectangles is equal to 
$\frac{b}{4}(4bc+2)=\frac{b(2bc+1)}{2}=\omega$.

Suppose now that $b\equiv 2 \pmod 4$.
We take $\mathcal{Q}_1=\varnothing$ if $b=6$ and
$$\mathcal{Q}_1=\left\{Q_i : i\in \left[0, \frac{(b-6)c}{4}-1\right] \right\}$$
if $b\geqslant 10$.
We have
$\E(\mathcal{Q}_1)=[1, (b-6)c]\cup 
[(b+6)c+1, 2bc]$.
For all $j\in [0,c-1]$, we define
$$S_j=\begin{array}{|c|c|c|c|} \hline
(b-6)c+4j+1 &  (b+6)c-4j-1 &
(b-6)c+4j+3 & (b+6)c-4j-3 \\\hline
(b+6)c-4j & (b-6)c+4j+2 &
(b+6)c-4j-2 &  (b-6)c+4j+4  \\ \hline
\end{array}.$$
We note that the sum of the elements of each column of $S_j$ is equal to $\delta=2bc+1$.
The sum of the elements of the first row of $S_j$ is $4bc$, while the sum of the elements of the second row is $4bc+4$.
Setting
$$\mathcal{S}=\left\{ S_j : j \in
[0,c-1]\right\},$$
we obtain
$\E(\mathcal{S})=[(b-6)c+1, 
(b-2)c] \cup
[(b+2)c+1,    (b+6)c ]$.
For all $h\geqslant 0$, we define
$$\begin{array}{rcl}
T_{2h} & = &\begin{array}{|c|c|} \hline
(b+2)c-4h &  (b-2)c+4h+3  \\\hline
(b-2)c+4h+1 & (b+2)c-4h-2 \\ \hline
\end{array},\\[12pt]
T_{2h+1} & = & \begin{array}{|c|c|} \hline
(b+2)c-4h-1 &  (b-2)c+4h+4  \\\hline
(b-2)c+4h+2 & (b+2)c-4h-3 \\ \hline
\end{array}.
\end{array}$$
Finally, we define
$$T' = \begin{array}{|c|c|} \hline
bc+2 & bc+1  \\\hline
bc-1 &  bc \\ \hline
\end{array}.$$
We notice that the sum of the elements of each column of $T_{2h},T_{2h+1},T'$ is equal to $\delta=2bc+1$.
The sum of the elements of the first row of each of these blocks is equal to 
$2bc+3$, while the sum of the elements of the second row is equal to $2bc-1$.
So, we take
$$\mathcal{T}=\left\{
\begin{array}{ll}
\left\{T_{2h},T_{2h+1}: h \in \left[0,\frac{c-2}{2}\right]\right\}&
\text{ if } c \text{ is even},\\[4pt]
\left\{T'\right\}\cup\left\{T_{2h},
T_{2h+1}: h \in \left[0,\frac{c-3}{2}\right]\right\}
& \text{ if } c \text{ is odd}.
\end{array}\right.$$
It follows that
$\E(\mathcal{T})= [(b-2)c+1,(b+2)c ]$.

To get an $\MRS(2,b;c)$ we take
$c$ rectangles obtained by juxtaposition, for each of them,  of $\frac{b-6}{4}$ distinct elements of $\mathcal{Q}_1$, one element of $\mathcal{S}$ and one element of  $\mathcal{T}$. In this way, the sum of the elements of the first row is  equal to $\frac{(4bc+2)(b-6)}{4}+4bc+(2bc+3)=
\frac{b(2bc+1)}{2}=\omega$
and the sum of the elements of the second row
is equal to $\frac{(4bc+2)(b-6)}{4}+(4bc+4)
+(2bc-1) =\omega$.
  \end{proof}

As well as Definition \ref{MRS} generalizes Definition \ref{defMRS}, 
we can also generalize Definition 
\ref{GMRS}.

\begin{defn}
Let $m,n,s,k,c$ be five positive integers.
A $\Gamma$-magic rectangle set $$\MRS_\Gamma(m,n;s,k;c)$$ is a set
of $c$ partially filled $m\times n$ arrays
with entries in an abelian group $\Gamma$ of order $nkc$ such that
    \begin{itemize}
        \item[{\rm (a)}] every element of $\Gamma$ appears exactly once and in a unique array;
        \item[{\rm (b)}] for every array, each row contains exactly $s$ filled cells and each column contains exactly $k$ filled cells;
        \item[{\rm (c)}] there exist two (not necessarily distinct) elements $\omega,\delta \in \Gamma$ such that, for every array, the sum of the entries of each row is equal to $\omega$ and the sum of the entries of each column is equal to $\delta$.
    \end{itemize}
\end{defn}

Our final goal is to provide necessary and sufficient conditions for the existence of an $\MRS_\Gamma(m,n;s,k;c)$ and this paper is a first attempt in this direction.

    
\begin{oss}\label{OssCic}
Taking a magic rectangle set $\MRS(m,n;s,k;c)$, 
it is possible to construct a
$\Z_{nkc}$-magic rectangle set
$\MRS_{\Z_{nkc}}(m,n;s,k;c)$ simply by applying the canonical epimorphism $\Z \to\Z_{nkc}$.
\end{oss}

\begin{oss}\label{tr}
An $\MRS_\Gamma(m,n;s,k;c)$ exists if and only if an $\MRS_\Gamma(n,m; k,s ;c)$ exists.
\end{oss}

\begin{lem}
    \label{eq:etichetta24}
     If the product $nc$ is even and  $k$ is odd, then
     there is no 
     $\MRS_\Gamma(m,n;s,k;c)$
     whenever
     $\Gamma$ is an abelian group of order $nkc$ having exactly one involution. 
\end{lem}

\begin{proof}
For the sake of contradiction, suppose there exists an $\MRS_\Gamma(m,n;s,k;c)$ for some abelian group $\Gamma$ of order $nkc$ having exactly one involution.  
We can construct a $\Gamma$-distance magic complete $nc$-partite graph
$G$ where each part has cardinality $k$
and labeling the vertices of each part
with the  entries of the columns of the elements of the
$\Gamma$-magic rectangle set (ignoring the empty cells).
The graph $G$ is $k(nc-1)$-regular, where $k(nc-1)$ is odd,
contradicting Theorem~\ref{th:graforegolare}.
\end{proof}



\begin{cor}\label{skodd1}
      If $s$ or $k$ is odd and  $\Gamma$ is an abelian group   of order $nkc$ having exactly one involution, then there is no $\MRS_\Gamma(m,n;s,k;c)$.
\end{cor}

\begin{proof}
 Since $\Gamma$ has an involution, 
 $|\Gamma|=nkc$ is even.
If $k$ is odd, then we apply Lemma 
    \ref{eq:etichetta24}.
If $s$ is odd, then we apply Remark \ref{tr}.    
\end{proof}

Easily adapting the proof of
\cite[Observation 3.8]{CH21b}, we have the following nonexistence result.

\begin{prop}\label{no_odd}
Let $\Gamma$ be an abelian group of order $2nc$. 
If $ns$ is odd, then there is no 
$\MRS_\Gamma(m,n; s,2;c)$. 
\end{prop}

\section{On  Cichacz-Hinc conjecture}

In this section,
we consider the validity of Conjecture
\ref{ConCH}.
First of all, 
Corollary \ref{skodd1} and Proposition \ref{no_odd} prove the necessity of the conditions.
Next, we consider the case when $c$ is odd.
This case was not explicitly considered by Cichacz and Hinc, but it can be obtained by adapting the proof of 
\cite[Theoren 3.11]{CH21b}.
To this purpose, we recall some auxiliary lemmas.

\begin{lem}\cite[Lemma 3.3]{CH21b}\label{lem:esistenzainsiemirettmagici}
    Let $\Gamma$ be an abelian group of order $abc_1c_2$ that can be decomposed as $ \Gamma\cong\Gamma_0\oplus \Phi$ for some group $\Phi\in\Upsilon$ of order $c_2$. 
    If there exists a $\Gamma_0$-magic rectangle set $\MRS_{\Gamma_0}(a,b;c_1)$, where the row sums are all equal to  $\omega_0$ and the column sums are all equal to  $\delta_0$, then there exists a $\Gamma$-magic rectangle set $\MRS_\Gamma(a,b;c_1c_2)$, where
    the row sums are all equal to $\omega=(\omega_0,0_\Phi)$ and the column sums are all equal to $\delta=(\delta_0,0_\Phi)$.
\end{lem}

\begin{lem}\cite[Lemma 3.4]{CH21b}\label{lem:2k+1}
    Let $\Gamma\cong\Lambda \oplus\Z_{2k+1}$ for some finite abelian group $\Lambda$. Let $h$ be a divisor of $2k+1$ and let $\Gamma_0\cong\Lambda\oplus\langle [h]_{2k+1}\rangle$. 
    If there exists an  $\MRS_{\Gamma_0}(a,b;c_1)$, where the row sums are all equal to $\omega_0$ and the column sums are all 
    equal to $\delta_0$, then there exists an 
    $\MRS_\Gamma(a,b;c_1h)$ where the row sums are all equal to $\omega_0+\left(0_\Lambda,\left[\frac{b(h+1)}{2}\right]_{2k+1}\right)$
    and the column sums are all equal to $\delta_0+\left(0_\Lambda,
    \left[\frac{a(h+1)}{2}\right]_{2k+1}\right)$.
\end{lem}

\begin{lem}\cite[Lemma 3.6]{CH21b}\label{lem:esistenza con primi dispari}
    If there exist two odd primes $p_1,p_2$ such that $p_1\mid a$ and $p_2\mid b$ and, furthermore, $\Gamma\in\Upsilon$ is an abelian group of order $abc$, then there exists an 
    $\MRS_\Gamma(a,b;c)$.
\end{lem}

\begin{lem}\cite[Lemma 3.7]{CH21b}\label{lem:p primo e più di 1 inv}
Let $p$ be a prime. There exists an 
$\MRS_\Gamma(p,2^\alpha;1)$ if and only if the abelian group $\Gamma$ of order $2^\alpha p$ has more than one involution.
\end{lem}

\begin{prop}
Let $\alpha,c$ be two integers with $\alpha\geqslant 2$ and $c \geqslant 1$. 
Let $\Gamma \in \Upsilon$ be an abelian group of order $(2\ell+1) 2^\alpha c$.
If $c$ is odd, then there exists an $\MRS_\Gamma(2\ell+1, 2^\alpha;c)$.
\end{prop}

\begin{proof}
   Since $\Gamma \in \Upsilon$ and  $|\Gamma|=(2\ell+1)2^\alpha c\equiv 0 \pmod 4$, the group $\Gamma$ contains more than one involution. Furthermore, since $c$ is odd, the product  $(2\ell+1)c$ is odd. Then, there exists an odd prime  divisor $p$ of $(2\ell+1)c$. 
   So, we can decompose $\Gamma$ as $\Gamma\cong\Z_{p^\beta}\oplus\Phi\oplus\Delta$, where $\beta\geqslant 1$, $|\Phi|=2^\alpha\geqslant 4$, $\Phi\in\Upsilon$, $|\Delta|=\frac{(2\ell+1)c}{p^\beta}\geqslant 1$ and $p\nmid |\Delta|$. The subgroup $\Psi=\langle [p^{\beta -1}]_{p^\beta} \rangle$ is cyclic of order $p$. 
   By Lemma \ref{lem:p primo e più di 1 inv} there exists an $\MRS_{\Psi\oplus\Phi}(p,2^\alpha;1)$. 
   Moreover, by Lemma \ref{lem:2k+1} there exists an $\MRS_{\Z_{p^\beta}\oplus\Phi}(p,2^\alpha;p^{\beta -1})$. Lemma \ref{lem:esistenzainsiemirettmagici} implies the existence of an  $\MRS_\Gamma\left(p,2^\alpha;\frac{(2\ell+1)c}{p}\right)$.
   Now, we take $\frac{2\ell+1}{p}$ distinct elements of $\MRS_\Gamma\left(p,2^\alpha;\frac{(2\ell+1)c}{p}\right)$, juxtaposing them to form a $(2\ell+1)\times 2^\alpha$ array.
   By repeating this process $c$ times, we get an $\MRS_\Gamma(2\ell+1,2^\alpha;c)$.
\end{proof}

We now prove the existence of an
$\MRS_\Gamma(2\ell + 1, 8; 4h+ 2)$ with $\Gamma \in \Upsilon$, providing in this way a further evidence for the validity of Conjecture \ref{ConCH}.
We start with a more general result.

\begin{lem}\label{2righe}
Let $\Omega$ be a subset of $\Z_{2a}$ such that if $g \in \Omega$ then $|g|>2$ and
$-g \in \Omega$.
Let $\Psi$ be an abelian group of order
$2^{\alpha}$ with $\alpha \geqslant 2$.
Then, there exists a set $\mathcal{T}$ consisting of $|\Omega|/2$ zero-sum blocks of size $2\times 2^{\alpha}$ such that 
$$\E(\mathcal{T})=\{(x,y): x \in \Omega,\; y \in \Psi\}
\subseteq \Z_{2a}\oplus \Psi.$$
\end{lem}

\begin{proof}
Let $2b\geqslant 2$ be the exponent of $\Psi$: then, we can write 
$\Psi=\Z_{2b}\oplus \Phi$, where $\Phi$ is a suitable abelian $2$-group of order $h=\frac{2^{\alpha}}{2b}$.
Let $(g_1,g_2,\ldots,g_{h})$ be any ordering of the elements of $\Phi$.
Given $x \in \Omega$ and $0\leqslant j\leqslant 2b-1$, 
let $A_{x,j}$ be the following $1\times h$ block:
$$A_{x,j}=\begin{array}{|c|c|c|c|}\hline
(x, [j]_{2b}, g_1) & (x, [j]_{2b}, g_2) & 
\cdots &  (x, [j]_{2b}, g_h) \\\hline
\end{array}.$$
The sum of the entries of this block is equal to 
$(hx, [hj]_{2b}, \sigma(\Phi))$, where $\sigma(\Phi)$ denotes the sum of the elements of $\Phi$. 

We recall that, by hypothesis, the elements of  $\Omega$ can be partitioned
into pairs of type $\{x,-x\}$ such that $x\ne -x$. 
So, fixing a such pair $\{x,-x\}$ of elements in $\Omega$, we build the block
$$T_x=\begin{array}{|c|c|c|c|c|c|c|c|}\hline
A_{x,0} & A_{x,1} & \cdots & A_{x,b-1} &
-A_{x,b} & -A_{x,b+1} & \cdots & -A_{x,2b-1}  \\\hline
-A_{x,0} & -A_{x,1} & \cdots & -A_{x,b-1} & A_{x,b} & A_{x,b+1} & \cdots & A_{x,2b-1}   \\ \hline
\end{array}$$
of size $2\times 2^\alpha$. 
Call $\mathcal{T}$ the set of these $|\Omega|/2$ blocks,
obtained by varying the pairs $\{x,-x\}$.
We observe that
$$\bigcup_{j=0}^{2b-1}\E(A_{x,j})=\{(x,y,z):
y \in\Z_{2b},\; z\in \Phi\}\subseteq \Z_{2a}\oplus \Z_{2b}\oplus \Phi$$
and that 
$\E(T_x)=\{(x,w),(-x,w): w \in \Psi\}\subseteq \Z_{2a}\oplus\Psi$.
Furthermore, the sum of the elements of each column of $T_x$ is $([0]_{2a},[0]_{2b},0_\Phi)$. 

Let $\varepsilon$ be the sum of the elements of the first row of $T_x$. Then, 
$$\begin{array}{rcl}
\varepsilon & =& \sum\limits_{j=0}^{b-1} \left(hx, [hj]_{2b}, \sigma(\Phi)\right)
-\sum\limits_{j=b}^{2b-1} \left(hx, [hj]_{2b}, \sigma(\Phi)\right)\\[4pt]
& = & \left( hbx, [h b(b-1)/2]_{2b} ,b \sigma(\Phi)\right) -
\left( hbx, [hb(3b-1)/2]_{2b}  , b\sigma(\Phi)\right) \\[4pt]
 & =& \left( [0]_{2a}, [-hb^2]_{2b} , 0_\Phi  \right).
\end{array}$$
Since $bh=2^{\alpha-1}\geqslant 2$ is an even integer,
we can write  $[-hb^2]_{2b}=[(2b)(-hb)/2]_{2b}=[0]_{2b}$.  We conclude that 
$\varepsilon =([0]_{2a},[0]_{2b},0_\Phi)$. 
Similarly, the sum of the elements of the second row is $-\varepsilon$ and so  $\mathcal{T}$ satisfies the requested properties.
\end{proof}

We now need a series of auxiliary lemmas.

\begin{lem}\label{lem:4r+4}
    Let $r\geqslant 3$ be an odd integer and let $\Gamma\cong\Z_{4r}\oplus\Z_4$. Then, there exists a zero-sum  $\MRS_\Gamma(r,8;2)$.
\end{lem}

\begin{proof}
First of all, we write the two elements
       $$\begin{array}{|c|c|} \hline
        S_{1,1} & S_{1,2} \\ \hline
        \end{array}\equad 
        \begin{array}{|c|c|} \hline
        S_{2,1} & S_{2,2} \\ \hline
        \end{array} $$
        of an $\MRS_{\Z_{12}\oplus\Z_4}(3,8;2)$, where
        $$S_{1,1}=\begin{array}{|c|c|c|c|cccc}\hline
 ( [0]_{12} , [0]_4) &  ([0]_{12} , [1]_4) & ([7]_{12} , [1]_4) & ([7]_{12} , [3]_4) \\\hline
 ( [2]_{12} , [1]_4) &  ([2]_{12} , [2]_4) & ([4]_{12} , [2]_4) & ([4]_{12} , [3]_4)\\\hline
 ([10]_{12} , [3]_4) & ([10]_{12} , [1]_4) & ([1]_{12} , [1]_4) & ([1]_{12} , [2]_4) \\\hline
        \end{array},$$
        $$S_{1,2}=\begin{array}{|c|c|c|c|}\hline
 ( [8]_{12} , [1]_4) & ( [8]_{12} , [2]_4) & ([9]_{12} , [3]_4) & ([9]_{12} , [1]_4)    \\\hline
 ( [6]_{12} , [3]_4) & ( [6]_{12} , [0]_4) & ([0]_{12} , [2]_4) & ([0]_{12} , [3]_4)  \\ \hline
 ([10]_{12} , [0]_4) & ([10]_{12} , [2]_4) & ([3]_{12} , [3]_4) & ([3]_{12} , [0]_4)\\\hline
        \end{array},$$
     $$S_{2,1}=\begin{array}{|c|c|c|c|}\hline
 ( [1]_{12} , [0]_4) & ([1]_{12} , [3]_4) & ([11]_{12} , [2]_4) & ([11]_{12} , [1]_4) \\\hline
 ( [2]_{12} , [0]_4) & ([2]_{12} , [3]_4) & ( [8]_{12} , [0]_4) & ( [8]_{12} , [3]_4) \\\hline
 ( [9]_{12} , [0]_4) & ([9]_{12} , [2]_4) & ( [5]_{12} , [2]_4) & ( [5]_{12} , [0]_4)\\\hline
        \end{array},$$
        $$S_{2,2}=\begin{array}{|c|c|c|c|}\hline
 ( [7]_{12} , [0]_4) & ( [7]_{12} , [2]_4) & ([5]_{12} , [1]_4) & ([5]_{12} , [3 ]_4) \\\hline
 ([11]_{12} , [3]_4) & ([11]_{12} , [0]_4) & ([3]_{12} , [2]_4) & ([3]_{12} , [1 ]_4) \\\hline
 ( [6]_{12} , [1]_4) & ( [6]_{12} , [2]_4) & ([4]_{12} , [1]_4) & ([4]_{12} , [0 ]_4) \\\hline
        \end{array}.$$       
Note that each row sum and each column sum is equal to
$([0]_{12},[0]_4)$.
    
    So, assume $r\geqslant 5$. We construct the first five rows of the two arrays of an $\MRS_{\Z_{4r}\oplus \Z_4}(r,8;2)$ by taking
        $$\begin{array}{|c|c|} \hline
        R_{1,1} & R_{1,2} \\ \hline
        \end{array}\equad 
        \begin{array}{|c|c|} \hline
        R_{2,1} & R_{2,2} \\ \hline
        \end{array}, $$
        where
        $$R_{1,1}=\begin{array}{|c|c|c|c|}\hline
        ([1]_{4r},[0]_4) & ([0]_{4r},[3]_4) & ([3r-1]_{4r},[0]_4) & ([r]_{4r},[1]_4)\\ \hline
        ([r-1]_{4r},[1]_4) & ([3r]_{4r},[3]_4) & ([1]_{4r},[3]_4) & ([0]_{4r},[1]_4)\\ \hline
        ([r]_{4r},[2]_4) & ([3r+2]_{4r},[3]_4) & ([0]_{4r},[2]_4) & ([2r-2]_{4r},[1]_4)\\ \hline
        ([2]_{4r},[3]_4) & ([2r-2]_{4r},[3]_4) & ([4r-1]_{4r},[1]_4) & ([1]_{4r},[1]_4)\\ \hline
        ([2r-2]_{4r},[2]_4) & ([0]_{4r},[0]_4) & ([r+1]_{4r},[2]_4) & ([r+1]_{4r},[0]_4)\\ \hline
        \end{array},$$
        $$R_{1,2}=\begin{array}{|c|c|c|c|}\hline
        ([4r-1]_{4r},[0]_4) & ([2r]_{4r},[3]_4) & ([3r+1]_{4r},[0]_4) & ([3r]_{4r},[1]_4)\\ \hline
        ([3r+1]_{4r},[1]_4) & ([r]_{4r},[3]_4) & ([2r-1]_{4r},[3]_4) & ([2r]_{4r},[1]_4)\\ \hline
        ([3r]_{4r},[2]_4) & ([r-2]_{4r},[3]_4) & ([2r]_{4r},[2]_4) & ([2]_{4r},[1]_4)\\ \hline
        ([4r-2]_{4r},[3]_4) & ([2r+2]_{4r},[3]_4) & ([2r+1]_{4r},[1]_4) & ([2r-1]_{4r},[1]_4)\\ \hline
        ([2r+2]_{4r},[2]_4) & ([2r]_{4r},[0]_4) & ([3r-1]_{4r},[2]_4) & ([r-1]_{4r},[0]_4)\\ \hline
        \end{array},$$
     $$R_{2,1}=\begin{array}{|c|c|c|c|}\hline
        ([2]_{4r},[0]_4) & ([3r+2]_{4r},[0]_4) & ([r-1]_{4r},[3]_4) & ([r+1]_{4r},[1]_4)\\ \hline
        ([r-1]_{4r},[2]_4) & ([r+2]_{4r},[0]_4) & ([3r-2]_{4r},[3]_4) & ([r+1]_{4r},[3]_4)\\ \hline
        ([1]_{4r},[2]_4) & ([2r-1]_{4r},[2]_4) & ([r+2]_{4r},[2]_4) & ([r-2]_{4r},[2]_4)\\ \hline
        ([r]_{4r},[0]_4) & ([2r-1]_{4r},[0]_4) & ([4r-1]_{4r},[3]_4) & ([4r-2]_{4r},[1]_4)\\ \hline
        ([2r-2]_{4r},[0]_4) & ([4r-2]_{4r},[2]_4) & ([3r+2]_{4r},[1]_4) & ([r+2]_{4r},[1]_4)\\ \hline
        \end{array},$$
        $$R_{2,2}=\begin{array}{|c|c|c|c|}\hline
        ([4r-2]_{4r},[0]_4) & ([r-2]_{4r},[0]_4) & ([3r+1]_{4r},[3]_4) & ([3r-1]_{4r},[1]_4)\\ \hline
        ([3r+1]_{4r},[2]_4) & ([3r-2]_{4r},[0]_4) & ([r+2]_{4r},[3]_4) & ([3r-1]_{4r},[3]_4)\\ \hline
        ([4r-1]_{4r},[2]_4) & ([2r+1]_{4r},[2]_4) & ([3r-2]_{4r},[2]_4) & ([3r+2]_{4r},[2]_4)\\ \hline
        ([3r]_{4r},[0]_4) & ([2r+1]_{4r},[0]_4) & ([2r+1]_{4r},[3]_4) & ([2r+2]_{4r},[1]_4)\\ \hline
        ([2r+2]_{4r},[0]_4) & ([2]_{4r},[2]_4) & ([3r-2]_{4r},[1]_4) & ([r-2]_{4r},[1]_4)\\ \hline
        \end{array}.$$
        Note that the elements appearing in these rows are exactly the elements of the set
        $$\{([a]_{4r}, x): a \in A,\;
        x\in  \Z_4\}\subseteq \Z_{4r}\oplus \Z_4,$$
        where $A=\{ 0,1, 2, r-2, r-1, r, r+1, r+2,
        2r-2, 2r-1, 2r, 2r+1,  2r+2,
        3r-2, 3r-1, 3r,3r+1, 3r+2, 4r-2, 4r-1\}$. 
        Furthermore, each row and column sum is equal to $([0]_{4r}, [0]_4)$.

        Next, we take the subset $\Omega=\Z_{4r}\setminus
        \{[a]_{4r}: a \in A\}$ of cardinality $4r-20$. By applying Lemma \ref{2righe}, we get a set $\mathcal{T}=\{T_i: 1\leqslant i\leqslant 2r-10\}$  of $2r-10$ zero-sum blocks of size $2\times 4$ such that
        $\E(\mathcal{T})=\{(x,y): x \in \Omega,\; y \in \Z_4\}$.

        Hence, we get two rows $\begin{array}{|c|c|}\hline
        T_{2j-1} & T_{2j} \\ \hline
        \end{array}$ by juxtaposing two distinct blocks of $\mathcal{T}$. 
        Finally, we complete the construction of the two arrays in  $\MRS_{\Z_{4r}\oplus\Z_4}(r,8;2)$
        by taking, for each one, $\frac{r-5}{2}$ of  these pairs of rows. In fact, each element of $\Z_{4r}\oplus\Z_4$ appears once and in a unique array and, in both arrays, the row and column sums are all equal to $\left([0]_{4r},[0]_4\right)$.
\end{proof}

\begin{lem}\label{lem:2r+8}
    Let $r\geqslant 3$ be an odd integer and let $\Gamma\cong\Z_{2r}\oplus\Z_8$. Then, there exists a zero-sum $\MRS_\Gamma(r,8;2)$.
\end{lem}

\begin{proof}
        First of all we write the two elements of an  $\MRS_{\Z_6\oplus\Z_8}(3,8;2)$:
        
        \begin{footnotesize}
        $$\begin{array}{|c|c|c|c|c|c|c|c|}\hline
        ([0]_6,[0]_8) & ([5]_6,[6]_8) & ([1]_6,[0]_8) & ([0]_6,[5]_8) & ([0]_6,[2]_8) & ([5]_6,[3]_8) & ([1]_6,[5]_8) & ([0]_6,[3]_8)\\ \hline
        ([1]_6,[6]_8) & ([5]_6,[4]_8) & ([3]_6,[6]_8) & ([3]_6,[4]_8) & ([1]_6,[1]_8) & ([5]_6,[1]_8) & ([3]_6,[0]_8) & ([3]_6,[2]_8)\\ \hline
        ([5]_6,[2]_8) & ([2]_6,[6]_8) & ([2]_6,[2]_8) & ([3]_6,[7]_8) & ([5]_6,[5]_8) & ([2]_6,[4]_8) & ([2]_6,[3]_8) & ([3]_6,[3]_8)\\ \hline
        \end{array},$$
        $$\begin{array}{|c|c|c|c|c|c|c|c|}\hline
        ([3]_6,[1]_8) & ([1]_6,[2]_8) & ([4]_6,[3]_8) & ([0]_6,[6]_8) & ([3]_6,[5]_8) & ([1]_6,[7]_8) & ([2]_6,[7]_8) & ([4]_6,[1]_8)\\ \hline
        ([4]_6,[7]_8) & ([4]_6,[2]_8) & ([2]_6,[1]_8) &([2]_6,[5]_8) & ([4]_6,[4]_8) & ([4]_6,[6]_8) & ([4]_6,[0]_8) & ([0]_6,[7]_8)\\ \hline
        ([5]_6,[0]_8) & ([1]_6,[4]_8) & ([0]_6,[4]_8) & ([4]_6,[5]_8) & ([5]_6,[7]_8) &  ([1]_6,[3]_8) & ([0]_6,[1]_8) & ([2]_6,[0]_8)\\ \hline
        \end{array}.$$
        \end{footnotesize}
        
       \noindent Note that each row sum and each column sum is 
equal to $([0]_6,[0]_8)$.
        
        Now, we assume $r\geqslant 5$. 
      We construct the first five rows 
       of the two arrays in the set $\MRS_{\Z_{2r}\oplus\Z_8}(r,8;2)$ 
       by taking
        $$\begin{array}{|c|c|} \hline
        R_{1,1} & R_{1,2} \\ \hline
        \end{array}\equad 
        \begin{array}{|c|c|} \hline
        R_{2,1} & R_{2,2} \\ \hline
        \end{array}, $$
        where
        $$R_{1,1}=\begin{array}{|c|c|c|c|}\hline
        ([1]_{2r},[0]_8) & ([0]_{2r},[0]_8) & ([r-1]_{2r},[1]_8) & ([r]_{2r},[1]_8)\\ \hline
        ([r-1]_{2r},[0]_8) & ([r]_{2r},[0]_8) & ([1]_{2r},[1]_8) & ([0]_{2r},[2]_8)\\ \hline
        ([r]_{2r},[2]_8) & ([r+2]_{2r},[0]_8) & ([0]_{2r},[3]_8) & ([2r-2]_{2r},[0]_8)\\ \hline
        ([2]_{2r},[5]_8) & ([2r-2]_{2r},[7]_8) & ([2r-1]_{2r},[0]_8) & ([1]_{2r},[6]_8)\\ \hline
        ([2r-2]_{2r},[1]_8) & ([0]_{2r},[1]_8) & ([r+1]_{2r},[3]_8) & ([r+1]_{2r},[7]_8)\\ \hline
        \end{array},$$
        $$R_{1,2}=\begin{array}{|c|c|c|c|}\hline
        ([1]_{2r},[4]_8) & ([0]_{2r},[4]_8) & ([r-1]_{2r},[7]_8) & ([r]_{2r},[7]_8)\\ \hline
        ([r-1]_{2r},[4]_8) & ([r]_{2r},[4]_8) & ([1]_{2r},[7]_8) & ([0]_{2r},[6]_8)\\ \hline
        ([r]_{2r},[6]_8) & ([r+2]_{2r},[4]_8) & ([0]_{2r},[5]_8) & ([2r-2]_{2r},[4]_8)\\ \hline
        ([2]_{2r},[0]_8) & ([2r-2]_{2r},[5]_8) & ([2r-1]_{2r},[7]_8) & ([1]_{2r},[2]_8)\\ \hline
        ([2r-2]_{2r},[2]_8) & ([0]_{2r},[7]_8) & ([r+1]_{2r},[6]_8) & ([r+1]_{2r},[5]_8)\\ \hline
        \end{array},$$
        $$R_{2,1}=\begin{array}{|c|c|c|c|}\hline
        ([2]_{2r},[2]_8) & ([r-2]_{2r},[1]_8) & ([2r-1]_{2r},[3]_8) & ([r+1]_{2r},[4]_8)\\ \hline
        ([r-1]_{2r},[3]_8) & ([r+2]_{2r},[5]_8) & ([r-2]_{2r},[2]_8) & ([r+1]_{2r},[0]_8)\\ \hline
        ([1]_{2r},[3]_8) & ([2r-1]_{2r},[5]_8) & ([r+2]_{2r},[6]_8) & ([r-2]_{2r},[6]_8)\\ \hline
        ([r]_{2r},[5]_8) & ([2r-1]_{2r},[6]_8) & ([r-1]_{2r},[2]_8) & ([2]_{2r},[1]_8)\\ \hline
        ([2r-2]_{2r},[3]_8) & ([2]_{2r},[7]_8) & ([r+2]_{2r},[3]_8) & ([r-2]_{2r},[5]_8)\\ \hline
        \end{array},$$
        $$R_{2,2}=\begin{array}{|c|c|c|c|}\hline
        ([2]_{2r},[4]_8) & ([r-2]_{2r},[0]_8) & ([2r-1]_{2r},[1]_8) & ([r+1]_{2r},[1]_8)\\ \hline
        ([r-1]_{2r},[6]_8) & ([r+2]_{2r},[7]_8) & ([r-2]_{2r},[7]_8) & ([r+1]_{2r},[2]_8)\\ \hline
        ([1]_{2r},[5]_8) & ([2r-1]_{2r},[2]_8) & ([r+2]_{2r},[2]_8) & ([r-2]_{2r},[3]_8)\\ \hline
        ([r]_{2r},[3]_8) & ([2r-1]_{2r},[4]_8) & ([r-1]_{2r},[5]_8) & ([2]_{2r},[6]_8)\\ \hline
        ([2r-2]_{2r},[6]_8) & ([2]_{2r},[3]_8) & ([r+2]_{2r},[1]_8) & ([r-2]_{2r},[4]_8)\\ \hline
        \end{array}.$$
        Note that the elements appearing in these rows are exactly the elements of the set
        $$\{([a]_{2r}, x): a \in A,\;
        x\in \Z_8\}\subseteq \Z_{2r}\oplus \Z_8,$$
        where $A=\{0,1,2, r-2,
        r-1, r, r+1, r+2, 2r-2, 
        2r-1\}$. 
        Furthermore, each row and column sum is equal to  $([0]_{2r},[0]_8)$.

        Next, we take the subset 
        $\Omega=\Z_{2r}\setminus
        \{[a]_{2r}: a \in A\}$ 
        of cardinality $2r-10$. By applying Lemma \ref{2righe}, we get a set $\mathcal{T}$ of $r-5$ zero-sum blocks of size $2\times 8$ such that
         $\E(\mathcal{T})=\{(x,y): x \in \Omega,\; y \in \Z_8\}$.
       Finally, we complete the construction of the two arrays in $\MRS_{\Z_{2r}\oplus\Z_8}(r,8;2)$
        by taking, for each one, $\frac{r-5}{2}$ of  these blocks.
        In fact, each element of $\Z_{2r}\oplus\Z_8$ appears once and in a unique array and, in both arrays, the row and column sums are all equal to $\left([0]_{2r},[0]_8\right)$.
\end{proof}

\begin{lem}\label{lem:2r+2+4}
  Let $r\geqslant 3$ be an odd integer and let $\Gamma\cong\Z_{2r}\oplus\Z_2\oplus\Z_4$. Then, there exists a zero-sum  $\MRS_\Gamma(r,8;2)$.
\end{lem}

\begin{proof}
We proceed as in the proof of Lemma \ref{lem:2r+8}.
So, we write the two elements
       $$\begin{array}{|c|c|} \hline
        S_{1,1} & S_{1,2} \\ \hline
        \end{array}\equad 
        \begin{array}{|c|c|} \hline
        S_{2,1} & S_{2,2} \\ \hline
        \end{array} $$
        of an  $\MRS_{\Z_6\oplus\Z_2\oplus \Z_4}(3,8;2)$, where
        $$S_{1,1}=\begin{array}{|c|c|c|c|}\hline
([0]_6 ,[0]_2 ,[3]_4) & ([0]_6 ,[0]_2 ,[0]_4) & ([5]_6 ,[1]_2 ,[1]_4) & ([5]_6 ,[1]_2 ,[2]_4) 
\\ \hline
([1]_6 ,[0]_2 ,[3]_4) & ([1]_6 ,[1]_2 ,[0]_4) & ([5]_6 ,[0]_2 ,[1]_4) & ([5]_6 ,[1]_2 ,[3]_4) 
\\\hline
([5]_6 ,[0]_2 ,[2]_4) & ([5]_6 ,[1]_2 ,[0]_4) & ([2]_6 ,[1]_2 ,[2]_4) & ([2]_6 ,[0]_2 ,[3]_4) \\\hline
        \end{array},$$
      $$S_{1,2}=\begin{array}{|c|c|c|c|}\hline
([1]_6 ,[1]_2 ,[3]_4) & ([1]_6 ,[1]_2 ,[2]_4) & ([3]_6 ,[0]_2 ,[0]_4) & ([3]_6 ,[0]_2 ,[1]_4) \\ \hline
([0]_6 ,[1]_2 ,[1]_4) & ([0]_6 ,[1]_2 ,[3]_4) & ([3]_6 ,[1]_2 ,[2]_4) & ([3]_6 ,[1]_2 ,[3]_4) \\ \hline
([5]_6 ,[0]_2 ,[0]_4) & ([5]_6 ,[0]_2 ,[3]_4) & ([0]_6 ,[1]_2 ,[2]_4) & ([0]_6 ,[1]_2 ,[0]_4) \\ \hline
        \end{array},$$
        $$S_{2,1}=\begin{array}{|c|c|c|c|}\hline
([2]_6 ,[0]_2 ,[0]_4) & ([2]_6 ,[0]_2 ,[1]_4) & 
([4]_6 ,[0]_2 ,[0]_4) & ([4]_6 ,[1]_2 ,[2]_4) \\ \hline
([3]_6 ,[0]_2 ,[2]_4) & ([3]_6 ,[0]_2 ,[3]_4) & 
([4]_6 ,[0]_2 ,[1]_4) & ([4]_6 ,[0]_2 ,[2]_4)  \\\hline
([1]_6 ,[0]_2 ,[2]_4) & ([1]_6 ,[0]_2 ,[0]_4) & 
([4]_6 ,[0]_2 ,[3]_4) & ([4]_6 ,[1]_2 ,[0]_4)  \\\hline
        \end{array},$$
   $$S_{2,2}=\begin{array}{|c|c|c|c|}\hline
([2]_6 ,[0]_2 ,[2]_4) & ([2]_6 ,[1 ]_2,[3]_4) &
([4]_6 ,[1]_2 ,[1]_4) & ([4]_6 ,[1 ]_2,[3]_4)   \\ \hline
([1]_6 ,[1]_2 ,[1]_4) & ([1]_6 ,[0]_2 ,[1]_4) & 
([0]_6 ,[0]_2 ,[2]_4) & ([2]_6 ,[1]_2 ,[0]_4)  \\\hline
([3]_6 ,[1]_2 ,[1]_4) & ([3]_6 ,[1]_2 ,[0]_4) & 
([2]_6 ,[1]_2 ,[1]_4) & ([0]_6 ,[0]_2 ,[1]_4) \\\hline
        \end{array}.$$

Next, we assume $r\geqslant 5$ and we construct the first five rows of the two arrays of an $\MRS_{\Z_{2r}\oplus\Z_2\oplus\Z_4}(r,8;2)$ by taking
     $$\begin{array}{|c|c|} \hline
        R_{1,1} & R_{1,2} \\ \hline
        \end{array}\equad 
        \begin{array}{|c|c|} \hline
        R_{2,1} & R_{2,2} \\ \hline
        \end{array}, $$
   where

    \begin{footnotesize}
    $$R_{1,1}=\begin{array}{|c|c|c|c|}\hline
        ([1]_{2r},[0]_2,[2]_4) & ([0]_{2r},[0]_2,[1]_4) & ([r-1]_{2r},[0]_2,[0]_4) & ([r]_{2r},[0]_2,[3]_4)\\ \hline
        ([r-1]_{2r},[0]_2,[2]_4) & ([r]_{2r},[0]_2,[1]_4) & ([1]_{2r},[1]_2,[3]_4) & ([0]_{2r},[1]_2,[1]_4)\\ \hline
        ([r]_{2r},[1]_2,[0]_4) & ([r+2]_{2r},[0]_2,[1]_4) & ([0]_{2r},[1]_2,[2]_4) & ([2r-2]_{2r},[0]_2,[0]_4)\\ \hline
        ([2]_{2r},[0]_2,[3]_4) & ([2r-2]_{2r},[0]_2,[1]_4) & ([2r-1]_{2r},[0]_2,[0]_4) & ([1]_{2r},[0]_2,[0]_4)\\ \hline
        ([2r-2]_{2r},[1]_2,[1]_4) & ([0]_{2r},[0]_2,[0]_4) & ([r+1]_{2r},[0]_2,[3]_4) & ([r+1]_{2r},[1]_2,[0]_4)\\ \hline
        \end{array},$$
        $$R_{1,2}=\begin{array}{|c|c|c|c|}\hline
        ([1]_{2r},[1]_2,[0]_4) & ([0]_{2r},[1]_2,[3]_4) & ([r-1]_{2r},[0]_2,[1]_4) & ([r]_{2r},[0]_2,[2]_4)\\ \hline
        ([r-1]_{2r},[0]_2,[3]_4) & ([r]_{2r},[0]_2,[0]_4) & ([1]_{2r},[0]_2,[3]_4) & ([0]_{2r},[0]_2,[3]_4)\\ \hline
        ([r]_{2r},[1]_2,[1]_4) & ([r+2]_{2r},[0]_2,[2]_4) & ([0]_{2r},[0]_2,[2]_4) & ([2r-2]_{2r},[1]_2,[0]_4)\\ \hline
        ([2]_{2r},[1]_2,[1]_2) & ([2r-2]_{2r},[0]_2,[3]_4) & ([2r-1]_{2r},[0]_2,[2]_4) & ([1]_{2r},[1]_2,[2]_4)\\ \hline
        ([2r-2]_{2r},[1]_2,[3]_4) & ([0]_{2r},[1]_2,[0]_4) & ([r+1]_{2r},[0]_2,[0]_4) & ([r+1]_{2r},[0]_2,[1]_4)\\ \hline
        \end{array},$$
        $$R_{2,1}=\begin{array}{|c|c|c|c|}\hline
        ([r+2]_{2r},[1]_2,[1]_4) & ([2r-2]_{2r},[0]_2,[2]_4) & ([2r-1]_{2r},[1]_2,[3]_4) & ([r+1]_{2r},[1]_2,[3]_4)\\ \hline
        ([r-1]_{2r},[1]_2,[1]_4) & ([r+2]_{2r},[0]_2,[3]_4) & ([r-2]_{2r},[0]_2,[0]_4) & ([r+1]_{2r},[0]_2,[2]_4)\\ \hline
        ([1]_{2r},[0]_2,[1]_4) & ([2r-1]_{2r},[1]_2,[0]_4) & ([r+2]_{2r},[0]_2,[0]_4) & ([r-2]_{2r},[1]_2,[3]_4)\\ \hline
        ([r]_{2r},[1]_2,[3]_4) & ([r-1]_{2r},[1]_2,[2]_4) & ([2r-1]_{2r},[1]_2,[1]_4) & ([2]_{2r},[1]_2,[0]_4)\\ \hline
        ([r-2]_{2r},[1]_2,[2]_4) & ([2]_{2r},[0]_2,[1]_4) & ([2]_{2r},[0]_2,[0]_4) & ([r-2]_{2r},[1]_2,[0]_4)\\ \hline
        \end{array},$$
        $$R_{2,2}=\begin{array}{|c|c|c|c|}\hline
        ([r+2]_{2r},[1]_2,[3]_4) & ([2r-2]_{2r},[1]_2,[2]_4) & ([2r-1]_{2r},[0]_2,[1]_4) & ([r+1]_{2r},[1]_2,[1]_4)\\ \hline
        ([r-1]_{2r},[1]_2,[3]_4) & ([r+2]_{2r},[1]_2,[0]_4) & ([r-2]_{2r},[0]_2,[1]_4) & ([r+1]_{2r},[1]_2,[2]_4)\\ \hline
        ([1]_{2r},[1]_2,[1]_4) & ([2r-1]_{2r},[0]_2,[3]_4) & ([r+2]_{2r},[1]_2,[2]_4) & ([r-2]_{2r},[0]_2,[2]_4)\\ \hline
        ([r]_{2r},[1]_2,[2]_4) & ([r-1]_{2r},[1]_2,[0]_4) & ([2r-1]_{2r},[1]_2,[2]_4) & ([2]_{2r},[1]_2,[2]_4)\\ \hline
        ([r-2]_{2r},[0]_2,[3]_4) & ([2]_{2r},[1]_2,[3]_4) & ([2]_{2r},[0]_2,[2]_4) & ([r-2]_{2r},[1]_2,[1]_4)\\ \hline
        \end{array}.$$
        \end{footnotesize}
        
        Note that the elements appearing in these rows are exactly the elements of
        $$\{([a]_{2r}, x, y): a \in A,\;
        x\in \Z_2, \; y\in\Z_4\}\subseteq \Z_{2r}\oplus \Z_2\oplus \Z_4,$$
        where $A=\{0,1,2, r-2,
        r-1, r, r+1, r+2, 2r-2, 
        2r-1\}$. 
        Furthermore, each row and column sum is equal to  $([0]_{2r},[0]_2,[0]_4)$.

        Now, we take the subset $\Omega=\Z_{2r}\setminus
        \{[a]_{2r}: a \in A\}$ of cardinality $2r-10$. 
        Applying Lemma  
        \ref{2righe}, we get a set $\mathcal{T}$  of
        $r-5$ zero-sum blocks  of size $2\times 8$ such that 
        $\E(\mathcal{T})=\{(x,y,z): x \in \Omega,\; y \in \Z_2,\;z\in \Z_4\}$.
        As in the previous lemma, we complete each array with $\frac{r-5}{2}$ of these blocks.
 \end{proof}

\begin{lem}\label{lem:2r+2+2+2}
    Let $r\geqslant 3$ be an odd integer and let $\Gamma\cong\Z_{2r}\oplus\Z_2\oplus\Z_2\oplus\Z_2$. Then, there exists a zero-sum  $\MRS_\Gamma(r,8;2)$.
\end{lem}

\begin{proof}
First of all, we construct the first three rows of each array in  $\MRS_{\Z_{2r}\oplus\Z_2\oplus \Z_2\oplus \Z_2}(r,8;2)$ 
        by taking
        $$\begin{array}{|c|c|c|c|} \hline
        R_{1,1} & R_{1,2} & R_{1,3} & R_{1,4} \\ \hline
        \end{array}\equad 
        \begin{array}{|c|c|c|c|} \hline
        R_{2,1} & R_{2,2} & R_{2,3} & R_{2,4} \\ \hline
        \end{array}, $$
        where
            $$R_{1,1}=\begin{array}{|c|c|}\hline
            ([0]_{2r}, [0]_2, [1]_2, [0]_2) & ( [0]_{2r}, [0]_2, [0]_2, [1]_2 ) \\\hline
            ([r]_{2r}, [0]_2, [0]_2, [0]_2) & ( [r]_{2r}, [1]_2, [0]_2, [0]_2) \\\hline
            ([r]_{2r}, [0]_2, [1]_2 ,[0]_2) & ( [r]_{2r}, [1]_2, [0 ]_2,[1]_2) \\\hline
            \end{array},$$
            $$R_{1,2}=\begin{array}{|c|c|}\hline
            ( [2r-2 ]_{2r},[0] ,[1], [0] ) & ([2r-2]_{2r}, [0]_2 ,[0]_2, [1]_2) \\ \hline
            ( [1]_{2r}, [0]_2, [0]_2 ,[0]_2) & ( [1]_{2r}, [1]_2, [0]_2, [0]_2 )\\ \hline
            ( [1]_{2r}, [0]_2, [1]_2, [0]_2 ) & ([1]_{2r}, [1]_2, [0]_2, [1]_2) \\ \hline
            \end{array},$$
            $$R_{1,3}=\begin{array}{|c|c|}\hline
            ([2]_{2r}, [0]_2, [1]_2, [0]_2) & ( [2]_{2r}, [1]_2, [0]_2, [0]_2 ) \\\hline
            ([2r-1]_{2r}, [0 ]_2,[0]_2, [0]_2 )& ( [2r-1]_{2r}, [0]_2, [1]_2, [1 ]_2) \\\hline
            ([2r-1]_{2r}, [0]_2, [1 ]_2,[0]_2 ) & ([2r-1]_{2r}, [1 ]_2,[1]_2, [1]_2) \\\hline
            \end{array},$$
            $$R_{1,4}=\begin{array}{|c|c|}\hline
            ([0]_{2r}, [0]_2, [0]_2, [0]_2 )& ( [0]_{2r}, [1]_2, [1]_2, [0]_2) \\ \hline
            ( [2r-2]_{2r}, [1]_2, [0 ]_2,[1]_2) & ( [2]_{2r}, [1]_2, [1 ]_2,[0]_2) \\ \hline
            ( [2]_{2r}, [1]_2, [0 ]_2,[1]_2) &( [2r-2]_{2r}, [0]_2, [0]_2, [0]_2) \\ \hline
            \end{array},$$
            $$R_{2,1}=\begin{array}{|c|c|}\hline
            ([0]_{2r} , [1]_2 , [1]_2 , [1]_2) & ([0]_{2r} , [1 ]_2, [0]_2 ,[ 0]_2)   \\\hline
            ([r]_{2r} , [1]_2 , [1]_2 , [0]_2) & ([r]_{2r} , [1]_2 , [1]_2 , [1]_2 ) \\\hline
            ([r]_{2r} , [0 ]_2, [0]_2 , [1]_2) & ([r]_{2r} , [0]_2 , [1]_2 , [1]_2 ) \\\hline
            \end{array},$$
            $$R_{2,2}=\begin{array}{|c|c|}\hline
            ([2r-2]_{2r} , [1]_2 , [0]_2 , [0]_2) & ([2r-2 ]_{2r}, [1]_2 , [1]_2 , [1]_2)   \\\hline
            ([1]_{2r} , [0 ]_2, [1 ]_2, [1]_2) & ([1]_{2r} , [1]_2 , [1]_2 , [0]_2 )   \\\hline
            ([1 ]_{2r}, [1 ]_2, [1]_2 , [1]_2) & ([1 ]_{2r}, [0]_2 , [0]_2 , [1 ]_2)  \\\hline
            \end{array},$$
            $$R_{2,3}=\begin{array}{|c|c|}\hline
            ([2 ]_{2r},   [0 ]_2, [0]_2 ,[ 1]_2 ),&  ([2]_{2r} ,     [1 ]_2, [1]_2 , [1  ]_2)   \\\hline
            ([2r-1]_{2r} , [1 ]_2, [0 ]_2, [1 ]_2) &  ([2r-1 ]_{2r}, [0 ]_2, [0 ]_2, [1  ]_2 ) \\\hline
            ([2r-1 ]_{2r}, [1 ]_2, [0]_2 ,[ 0  ]_2) &  ([2r-1 ]_{2r}, [1]_2 , [1 ]_2, [0 ]_2  ) \\\hline
            \end{array},$$
            $$R_{2,4}=\begin{array}{|c|c|}\hline
            ([0]_{2r} ,   [0 ]_2, [1 ]_2, [1 ]_2)& ([0 ]_{2r},    [1]_2 , [0]_2 , [1  ]_2) \\\hline
            ([2r-2]_{2r} , [0]_2 , [1]_2 , [1]_2) & ([2  ]_{2r}  , [0 ]_2, [1]_2 , [1 ]_2 ) \\\hline
            ([2 ]_{2r}, [0 ]_2,   [ 0]_2 , [0]_2 ) & ([2r-2 ]_{2r}, [1 ]_2,[ 1 ]_2, [0]_2)   \\\hline
        \end{array}.$$
        Note that the elements appearing in these rows are exactly the elements of the set
        $$\{([a]_{2r}, x, y,z): a \in A,\;
        x,y,z\in \Z_2\}\subseteq \Z_{2r}\oplus \Z_2\oplus \Z_2\oplus \Z_2,$$
        where $A=\{0,1,2, r, 2r-2, 2r-1 \}$. 
        Furthermore, each row and column sum is equal to  $([0]_{2r},[0]_2,[0]_2,[0]_2)$.

        Next, we take the subset $\Omega=\Z_{2r}\setminus
        \{[a]_{2r}: a \in A\}$ of cardinality $2r-6$. Applying Lemma
        \ref{2righe}, we get a set $\mathcal{T}$ of $r-3$ zero-sum
        blocks of size  $2\times 8$ such that
        $\E(\mathcal{T})=\{(x,y,z,w): x \in \Omega,\; y \in \Z_2,\;z\in \Z_2,\; w\in  \Z_2\}$.
        We complete each array in  $\MRS_{\Z_{2r}\oplus\Z_2\oplus \Z_2\oplus \Z_2}(r,8;2)$ with $\frac{r-3}{2}$ of these blocks.
        In fact, each element of $\Z_{2r}\oplus\Z_2\oplus \Z_2\oplus \Z_2$ appears once and in a unique array and, in both arrays, each row and column sum is $\left([0]_{2r},[0]_2,[0]_2,[0]_2\right)$.
\end{proof}

    \begin{proof}[Proof of Theorem \ref{main}]
        Since $|\Gamma|=16(2\ell+1)(2h +1)$, the group $\Gamma$ can have one, three, seven or fifteen involutions. If $\Gamma$ has a unique involution, by Lemma \ref{eq:etichetta24} there is no $\MRS_{\Gamma}(2\ell+1,8;4h +2)$. 
        So, the $2$-component $\Delta$ of $\Gamma$ is
        $\Delta\cong\Z_4\oplus\Z_4$, 
        or $\Delta\cong\Z_2\oplus\Z_8$, or $\Delta\cong\Z_2\oplus\Z_2\oplus\Z_4$, or $\Delta\cong \Z_2\oplus\Z_2\oplus\Z_2\oplus\Z_2$. By the structure theorem of finite abelian groups, the group $\Gamma$ can be written as $$\Gamma\cong\Delta\oplus\Z_{q_1^{\alpha_1}}\oplus\ldots\oplus\Z_{q_w^{\alpha_w}},$$
        where $q_1, \ldots,q_w$ are (not necessarily distinct) odd primes and $(2\ell+1)(2h+1)=q_1^{\alpha_1}\cdots q_w^{\alpha_w}$.
        Without loss of generality, we can assume $2\ell+1=q_1\cdot q$ for some odd integer $q\geqslant 1$.

        As in the proof of \cite[Theorem 2.5]{CH21}, we take an $\MRS_{\Delta\oplus\Z_{q_1}}(q_1,8;2)$ as starting point and we transform this set into an  $\MRS_{\Gamma}(2\ell+1,8;4h +2)$ by applying Lemmas \ref{lem:esistenzainsiemirettmagici} and \ref{lem:2k+1}. More precisely, by Lemmas \ref{lem:4r+4}, \ref{lem:2r+8}, \ref{lem:2r+2+4} and  \ref{lem:2r+2+2+2}, there exists an $\MRS_{\Z_{q_1}\oplus \Delta}(q_1,8;2)$.
        Now, applying Lemma \ref{lem:2k+1}, we get the existence of an $\MRS_{\Delta\oplus\Z_{q_1^{\alpha_1}}}\left(q_1,8;2q_1^{\alpha_1-1}\right)$ (it can happen that $\alpha_1=1)$. Lemma \ref{lem:esistenzainsiemirettmagici} implies the existence of an  $\MRS_{\Gamma}\left(q_1,8;\frac{2(2\ell+1)(2h +1)}{q_1}\right)$. To construct each array in  $\MRS_{\Gamma}(2\ell+1,8;4h +2)$, we juxtapose $\frac{2\ell+1}{q_1}$  elements of $\MRS_{\Gamma}\left(q_1,8;\frac{2(2\ell+1)(2h +1)}{q_1}\right)$ in such a way to get a  $(2\ell+1)\times 8$ array.
    \end{proof}

\section{Futher constructions}

In this section we present some constructions of  $\MRS_\Gamma(m,n;s,k;c)$.
We recall that the conditions \eqref{nec}
are necessary for the existence of these objects.
Setting $d=\gcd(s,k)$, we can write
\begin{equation}\label{s1}
s=ds_1,\quad k=dk_1, \quad m=ek_1\equad n=es_1
\end{equation}
for a suitable integer $e\geqslant d$.

It will be convenient to denote an $\MRS_\Gamma(n,n; k,k;c)$ by $\MRS_\Gamma(n;k;c)$.

\begin{oss}
Keeping the previous notation, 
if there exists an $\MRS_\Gamma(e;d;s_1k_1c)$, then there exists an $\MRS_\Gamma(m,n;s,k;c)$.
\end{oss}

Given a square partially filled array  $A=(a_{i,j})$
of size $n$,
we define the diagonal $D_\ell$, 
where $0\leqslant \ell \leqslant n-1$, as the set of the cells
$$\{(i,j): j-i\equiv \ell \pmod n\}$$ of $A$. 
We say that an $\MRS_\Gamma(n;k;c)$ is diagonal if, for every array of this set, the filled cells are exactly those of  $k$ consecutive diagonals
(working modulo $n$ on the indices of $D_\ell$). 

\begin{ex}
Here an example of a diagonal $\MRS_{\Z_6\oplus \Z_2\oplus \Z_4}(6;4;2)$, where the entries  $xyz$ must be read as $([x]_6,[y]_2,[z]_4)$:
 $$\begin{array}{|c|c|c|c|c|c|}\hline
000 & 500 &     &     & 410 & 310 \\\hline
311 & 001 & 503 &     &     & 413 \\\hline
010 & 110 & 002 & 502 &     &      \\\hline
    & 013 & 111 & 003 & 501 &      \\\hline
    &     & 012 & 112 & 200 & 300  \\\hline
303 &     &     & 011 & 113 & 201  \\\hline
   \end{array},\quad
\begin{array}{|c|c|c|c|c|c|}\hline
202 & 302 &     &     & 212 & 512 \\\hline
513 & 203 & 301 &     &     & 211 \\\hline
412 & 312 & 400 & 100 &     &      \\\hline
    & 411 & 313 & 401 & 103 &     \\\hline
    &     & 210 & 510 & 402 & 102 \\\hline
101 &     &     & 213 & 511 & 403 \\\hline
   \end{array}.$$
\end{ex}

The importance of building a diagonal $\MRS(n;k;c)$ is justified by the following result.

\begin{thm}\label{sq->rt}\cite{MP3,MP4}
Let $m,n,s,k,c$ be five positive integers such that $2\leqslant s \leqslant n$, $2\leqslant k\leqslant m$ e $ms=nk$.
Let $\Gamma$ be an abelian group of order $nkc$. 
Let $d=\gcd(s,k)$.
If there exists a diagonal $\MRS_\Gamma\left( \frac{nk}{d}; d; c \right)$, then there exists an $\MRS_\Gamma ( m,n; s,k;c)$.
\end{thm}

More in general, given a  partially filled array $A=(a_{i,j})$
of size $m\times n$,
we set $e=\gcd(m,n)$ and define
the diagonal $D_\ell$, where
$0\leqslant \ell \leqslant e-1$, as the set of cells
$$\{(i,j): j-i\equiv \ell \pmod{e} \}$$
of $A$.
For instance, given
$$A=\begin{pmatrix}
1 & 2 & 3 & 4 & 5 & 6 & 7 & 8 & 9\\
11 & 12 & 13 & 14 & 15 & 16 & 17 & 18 & 19\\
21 & 22 & 23 & 24 & 25 & 26 & 27 & 28 & 29  
\end{pmatrix},$$
the elements of the diagonal $D_1$ are
$(2,13,24, 5,16,27, 8,19,21)$.

Fixed the row index $i$, there are exactly  $\frac{n}{e}$ integers $j$ in $[1,n]$ such that
$j\equiv i + \ell \pmod{e}$.
This means that every row of  $A$ contains exactly
$\frac{n}{e}$ cells belonging to the  diagonal $D_\ell$.
Similarly, fixed the column index $j$,
there exist exactly  $\frac{m}{e}$ integers $i$ in $[1,m]$ such that
$i\equiv j-\ell \pmod{e}$.
It follows that each column of $A$ contains exactly 
$\frac{m}{e}$ cells belonging to the diagonal $D_\ell$.
In particular, each diagonal $D_\ell$ contains
$\frac{n}{e} m=\frac{m}{e}n=\lcm(m,n)$ distinct elements.

\begin{oss}\label{buchi}
Given four integers $m,n,s,k$ such that
$1\leqslant s \leqslant n$, $1\leqslant k\leqslant m$ and $ms=nk$, we can take $d=\gcd(s,k)$ diagonals in an $m\times n$ array: in this way we have selected $s$ cells in each row and $k$ cells in each column.
\end{oss}

We would like to determine necessary and sufficient conditions for the existence of a diagonal $\MRS_\Gamma(n;k;c)$.

\begin{thm}\label{MRSdia}
Let $n,c$ be two integers such that $n\geqslant 2$ and $c\geqslant 1$. A diagonal $\MRS_{\Gamma}(n;2;c)$ exists if and only if  $\Gamma$ contains an element of order $n$.
\end{thm}

\begin{proof}
Suppose there exists a diagonal $\MRS_{\Gamma}(n;2;c)$.
Calling  $R_1,\ldots,R_c$ its elements,
it must exist  $\alpha,\omega \in \Gamma$ and suitable elements 
$z_1,\ldots,z_c \in \Gamma$ such that
$$R_i=\begin{array}{|c|c|c|c|c|}\hline
z_i & \omega-z_i & & & \\ \hline
     & \alpha+z_i & \omega -(\alpha+z_i) & & \\ \hline
     & & 2\alpha+z_i & \omega -(2\alpha+z_i) & \\ \hline
     & & & \ddots & \ddots \\ \hline
    \omega -((n-1)\alpha+z_i) & & & & (n-1)\alpha+z_i \\ \hline
    \end{array},$$
    where the sum of the entries of each row is equal to
$\omega$, while the sum of the entries of each column is equal to  $\omega+\alpha$. Considering the first column, we get the condition $n\alpha=0_\Gamma$. Since the elements
$z_i,\alpha+z_i,\ldots, (n-1)\alpha+z_i$ are pairwise distinct, we get  $|\alpha|=n$.
In particular, calling  $H$ the subgroup of  $\Gamma$ generated by $\alpha$, the elements of the diagonals  $D_0,D_1$ of $R_i$ are the elements of the cosets $H+z_i$ and
$H+\omega-z_i$.

Hence, if there exists a diagonal $\MRS_{\Gamma}(n;2;c)$, then $\Gamma$ must contain an element of order  $n$. We now show that this condition is also sufficient.

So, let $H$ be a cyclic subgroup of $\Gamma$ of order $n$ and let $\pi: \Gamma \to \frac{\Gamma}{H}$ be the canonical epimorphism. Since $\Gamma$ is an abelian group of order $2nc$, the quotient group $\frac{\Gamma}{H}$ is abelian of order $2c$. Hence, there exists an isomorphism  $\psi: \frac{\Gamma}{H} \to \Psi=\Z_{2a}\oplus \Z_{b_1}\oplus \ldots \oplus\Z_{b_k}$ for some $k\geqslant 0$,
$a\geqslant 1$ e $b_1,b_2,\ldots,b_k> 1$ such that  $c = a b_1\cdots b_k$. It follows that $\psi \circ \pi: \Gamma \to \Psi$ is an epimorphism. Let $\xi = ([-1]_{2a},[-1]_{b_1},\ldots,[-1]_{b_k}) \in \Psi$.
Write  $\Phi=\Z_{b_1}\oplus \ldots \oplus \Z_{b_k}$, so that $\Psi=\Z_{2a}\oplus \Phi$; furthermore, we set  $h=|\Phi|\geqslant 1$. Let $(g_1,\ldots,g_h)$ be any ordering of the elements of $\Phi$.
For any  $i=1,\ldots,a$ and any $j=1,\ldots,h$, let
$$y_{i,j}=( [i-1]_{2a}, g_j ) \in \Z_{2a}\oplus \Phi.$$

Given $1\leqslant i,i'\leqslant a$ and $1\leqslant j,j'\leqslant h$, suppose that  $y_{i,j}=\xi- y_{i',j'}$.
Then, $[i-1]_{2a}=[-1-i'+1]_{2a}$ implies $i+i' \equiv  1\pmod{2a}$,  in contradiction with $2\leqslant i+i'\leqslant 2a$. So, the set $\{y_{i,j}, \xi-y_{i,j}: 1 \leqslant i\leqslant a,\; 1\leqslant j\leqslant h\}$ consists of $2ah=2c$ pairwise distinct elements of $\Psi$.

Now, for any $ 1 \leqslant i\leqslant a$ and any $ 1\leqslant j\leqslant h$, we take an element
$x_{i,j}\in \Gamma$ such that $(\psi\circ\pi)(x_{i,j})=y_{i,j}$.
Furthermore, let $\omega \in \Gamma$ be such that $(\psi\circ\pi)(\omega)=\xi$.
It follows that the cosets  $H+x_{i,j}$, $H+\omega-x_{i',j'}$, where $i,i'\in \{1,\ldots, a\}$ and
 $ j,j'\in \{1,\ldots, h\}$, are pairwise disjoint.
In fact, this holds if and only if the  elements $\pi(x_{i,j}), \pi(\omega-x_{i',j'})$ are pairwise distinct in $\frac{\Gamma}{H}$. On the other hand, this holds, as previously seen, since the elements $\psi(\pi(x_{i,j}))=y_{i,j}$, $\psi(\pi(\omega-x_{i',j'}))=\xi-y_{i',j'}$ are pairwise distinct in $\Psi$.

Set $z_{h(i-1) + j} =x_{i,j}$ where $ 1 \leqslant i\leqslant a$ and $ 1\leqslant j\leqslant h$.
We can now use the elements $z_1,\ldots,z_c$
 to build the arrays $R_i$ as described at the beginning of this proof, that is, where $\E(R_i)=(H+z_i) \cup (H+\omega-z_i)$. These arrays $R_1,\ldots,R_c$ are  the elements of a diagonal $\MRS_\Gamma(n;2;c)$.
\end{proof}

\begin{ex}
We construct a diagonal $\MRS_{\Z_4\oplus\Z_8}(4;2;4)$. Following the proof of the previous theorem, we take $\alpha=([1]_4,[2]_8)$ and we denote by $H$ the subgroup of $\Z_4\oplus\Z_8$  generated by this element of order $4$. Then, $H=\{([1]_4,[2]_8), ([2]_4,[4]_8), ([3]_4,[6]_8), ([0]_4,[0]_8)\}$. 

    Let $\pi:\Z_4\oplus\Z_8\to\frac{\Z_4\oplus\Z_8}{H}$ be the canonical epimorphism such that $\pi\left(([x]_4,[y]_8)\right)=H+([x]_4,[y]_8)$. The group $\frac{\Z_4\oplus\Z_8}{H}$ is cyclic of order $8$, being generated, for example, by $H+([0]_4,[1]_8)$.
    Hence, there exists an isomorphism $\psi:\frac{\Z_4\oplus\Z_8}{H}\to\Z_8$   such that $\psi(H+([x]_4,[y]_8))=
    [2x-y]_8$. 
    Next, consider the epimorphism $\psi\circ\pi:\Z_4\oplus\Z_8\to\Z_8$.
    Let  $\xi=[-1]_8$. So, we can take $\omega=([0]_4,[1]_8)$. Moreover, we get the following elements $y_{i,1}\in\Z_8$:
    $$y_{1,1}=[0]_8, \; \; y_{2,1}=[1]_8, \; \; y_{3,1}=[2]_8, \; \; y_{4,1}=[3]_8.$$
    By computing the preimages of the elements $y_{i,1}$ via the epimorphism $\psi\circ\pi$, we obtain the elements $z_{\ell}$:
    $$z_1=([0]_4,[0]_8), \; \; z_2=([1]_4,[1]_8), \; \; z_3=([2]_4,[2]_8), \; \; z_4=([3]_4,[3]_8).$$
    Following the proof of the previous theorem, we can construct the arrays $R_1,\ldots,R_4$ of a diagonal  $\MRS_{\Z_4\oplus\Z_8}(4;2;4)$:
    $$R_1=\begin{array}{|c|c|c|c|}\hline
    ([0]_4,[0]_8) & ([0]_4,[1]_8) & & \\ \hline
     & ([1]_4,[2]_8) & ([3]_4,[7]_8) & \\ \hline
     & & ([2]_4,[4]_8) & ([2]_4,[5]_8) \\ \hline
    ([1]_4,[3]_8) & & & ([3]_4,[6]_8) \\ \hline
    \end{array},$$
    $$R_2=\begin{array}{|c|c|c|c|}\hline
    ([1]_4,[1]_8) & ([3]_4,[0]_8) & & \\ \hline
     & ([2]_4,[3]_8) & ([2]_4,[6]_8) & \\ \hline
     & & ([3]_4,[5]_8) & ([1]_4,[4]_8) \\ \hline
    ([0]_4,[2]_8) & & & ([0]_4,[7]_8) \\ \hline
    \end{array},$$
    $$R_3=\begin{array}{|c|c|c|c|}\hline
    ([2]_4,[2]_8) & ([2]_4,[7]_8) & & \\ \hline
     & ([3]_4,[4]_8) & ([1]_4,[5]_8) & \\ \hline
     & & ([0]_4,[6]_8) & ([0]_4,[3]_8) \\ \hline
    ([3]_4,[1]_8) & & & ([1]_4,[0]_8) \\ \hline
    \end{array},$$
    $$R_4=\begin{array}{|c|c|c|c|}\hline
    ([3]_4,[3]_8) & ([1]_4,[6]_8) & & \\ \hline
     & ([0]_4,[5]_8) & ([0]_4,[4]_8) & \\ \hline
     & & ([1]_4,[7]_8) & ([3]_4,[2]_8) \\ \hline
    ([2]_4,[0]_8) & & & ([2]_4,[1]_8) \\ \hline
    \end{array}.$$
    In fact, each element of $\Z_4\oplus\Z_8$ appears exactly once and in a unique array. Moreover, the sum of the entries in each row is equal to $([0]_4,[1]_8)=\omega$ and the sum of the entries in each column is equal to $([1]_4,[3]_8)=\omega +\alpha$. 
\end{ex}

\begin{cor}\label{MRS2b}
Let $b,n,c$ be three positive integers such that $2\leqslant 2b\leqslant n$. Let $\Gamma$ be an abelian group of order $2nbc$.
If $\Gamma$ contains an element of order $n$, then there exists a diagonal $\MRS_\Gamma(n;2b;c)$.
In particular, this holds if $n$ divides $\exp(\Gamma)$.
\end{cor}

\begin{proof}
Take a diagonal $\MRS_\Gamma(n;2;bc)$ whose existence follows from Theorem \ref{MRSdia}.
Let $R_1,R_2,$ $\ldots, R_{bc}$ be its elements and denote by
$\omega$ and $\delta$, respectively, the sums of the elements in each row and column of $R_i$.

Fixed an integer $\ell$, with $1\leqslant \ell \leqslant c$,
we construct an array $S_\ell$ whose diagonals
$D_{2j},D_{2j+1}$ are filled, respectively and respecting the ordering, with  the elements of the diagonals
$D_0,D_1$ of $R_{(\ell-1)b+j+1}$, for every $0\leqslant j \leqslant b-1$.

Calling $\mathcal{S}$ the set of the arrays
$S_1,\ldots,S_c$ so constructed, every element of  $\Gamma$ appears once and in a unique array of $\mathcal{S}$. Furthermore, the sum of the entries of each row of $S_\ell$ is equal to $b\omega$ and the sum of the entries of each column is equal to  $b\delta$.
This proves that $\mathcal{S}$ is a diagonal $\MRS_\Gamma(n;2b;c)$.

The last part of the statement follows from the fact that a finite abelian group always has an element $g$ such that $|g|=\exp(\Gamma)$.
\end{proof}

The condition that  $\Gamma$ contains an element of order $n$ is, in general, sufficient but not necessary for the existence of a diagonal $\MRS_\Gamma(n;k;c)$ with $k\geqslant 4$ even. In fact, we can always construct a such set when  $k\equiv 0 \pmod 4$.

\begin{prop}\label{diag4}
Let $b,n,c$ be three positive integers such that
$4\leqslant  4b\leqslant n$.
For every abelian group $\Gamma$ of order $4nbc$ there exists a diagonal $\MRS_\Gamma(n;4b;c)$.
\end{prop}

\begin{proof}
By Theorem \ref{MRSdia}, there exists an $\MRS_\Gamma(2,2; nbc)$: let
$R_1,\ldots,R_{nbc}$ be its elements and write
$$R_i=\begin{array}{|c|c|} \hline
x_i & y_i \\ \hline
z_i & w_i \\ \hline
\end{array}.$$
Let $\omega$ and $\delta$ be, respectively, the sum of the elements of each row and column. For all  $1\leqslant i\leqslant nbc$, consider the following $3\times 2$ partially filled array with entries
in $\Gamma$:
    $$B_i=\begin{array}{|c|c|}\hline
    x_i & y_i \\ \hline
     & \\ \hline
    z_i & w_i\\ \hline
    \end{array}.$$
    Define $\mathcal{B}=\left(B_1,\ldots, B_{nbc}\right)$. 
    
    We use these $3\times 2$ blocks to build partially filled arrays, as we now show. Take an empty array $A_1$ of size  $n\times n$; arrange in this array the first $n$ blocks of  $\mathcal{B}$ in such a way that the element $x_i$ of the cell $(1,1)$ of $B_i$ fill the cell  $(i,i)$ of $A_1$ (we work modulo $n$ on row and column indices). 
    In this way, we fill the diagonal  $D_{n-2},D_{n-1},D_0,D_1$. In particular, each row and each column has four filled cells. Looking at the rows, the elements belonging to the diagonals $D_0, D_1$ have sum equal to $\omega$ and the same holds for the elements belonging to the diagonals $D_{n-2},D_{n-1}$. 
    Looking at the columns, the elements belonging to the diagonals $D_0, D_{n-2}$ have sum equal to  $\delta$ and the same holds for the elements belonging to the diagonals $D_1, D_{n-1}$.  
    It follows that  $A_1$ has row and column sums equal, respectively, to $2\omega$ and $2\delta$.

    Applying this process $b$ times  (working with $B_{n+1},B_{n+2},\ldots,B_{2n}$ on the diagonals $D_2,D_3,$ $D_4,D_5$, and so on), we obtain a partially filled array $A_1$,  with $4b$ filled cells in each row and column and where the row and column sums are, respectively, equal to
    $2b\omega$ and $2b\delta$.
    Finally, we repeat this entire process $c$ times, 
    obtaining a set $\{A_1,\ldots,A_c\}$ of partially filled arrays.
    Furthermore, by construction of the set $\mathcal{B}$, 
    it is clear that  $\displaystyle\bigcup_{\ell =1}^{c}\E(A_\ell)=
    \E(\mathcal{B})$ is equal to $\Gamma$. 
    So, we have proved that the set $\{A_1,\ldots,A_c\}$ is a diagonal $\MRS_\Gamma(n;4b;c)$.
    \end{proof}

\begin{ex}
We construct a diagonal $\MRS_{\Z_4\oplus\Z_4\oplus\Z_2}(8;4;1)$.
    Given $\Gamma=\Z_4\oplus\Z_4\oplus\Z_2$, by Theorem \ref{MRSdia} there exists an $\MRS_{\Gamma}(2,2;8)$.
    Indeed, take the involution $\alpha=([0]_4,[0]_4,[1]_2)$ and consider the epimorphism $\varphi: \Gamma\to \Z_4 \oplus \Z_4$
    given by $\varphi(([x]_4,[y]_4,[z]_2))=([x]_4,[y]_4)$.
    Following the proof of the previous proposition, consider the following partially filled arrays $B_1,\ldots,B_8$, whose entries $xyz$ must be read as $([x]_4,[y]_4,[z]_2)$:
    $$B_1=\begin{array}{|c|c|}\hline
    000 & 330 \\ \hline
     & \\ \hline
    331 & 001 \\ \hline
    \end{array}, \; 
    B_2=\begin{array}{|c|c|}\hline
    010 & 320 \\ \hline
     & \\ \hline
    321 & 011 \\ \hline
    \end{array}, \; 
    B_3=\begin{array}{|c|c|}\hline
    020 & 310 \\ \hline
     & \\ \hline
    311 & 021 \\ \hline
    \end{array}, \; 
    B_4=\begin{array}{|c|c|}\hline
    030 & 300 \\ \hline
     & \\ \hline
    301 & 031 \\ \hline
    \end{array}, $$
    $$ B_5=\begin{array}{|c|c|}\hline
    100 & 230 \\ \hline
     & \\ \hline
    231 & 101 \\ \hline
    \end{array}, \; 
    B_6=\begin{array}{|c|c|}\hline
    110 & 220 \\ \hline
     & \\ \hline
    221 & 111 \\ \hline
    \end{array}, \;
    B_7=\begin{array}{|c|c|}\hline
    120 & 210 \\ \hline
     & \\ \hline
    211 & 121 \\ \hline
    \end{array}, \;
    B_8=\begin{array}{|c|c|}\hline
    130 & 200 \\ \hline
     & \\ \hline
    201 & 131 \\ \hline
    \end{array}.$$  
    Observe that the sum of the entries in each nonempty row is equal to $([3]_4,[3]_4,[0]_2)$ and the sum of 
    the entries in each column is equal to $([3]_4,[3]_4,[1]_2)$.

    We can now construct a diagonal $\MRS_{\Z_4\oplus\Z_4\oplus\Z_2}(8;4;1)$, as follows:
    $$\begin{array}{|c|c|c|c|c|c|c|c|}\hline
    000 &   330 &     &    &    &     & 211   & 121 \\ \hline
    131 & 010   & 320 &    &    &     &      & 201 \\ \hline
    331 & 001   & 020 & 310 &    &     &      &    \\ \hline
        & 321   & 011 & 030 & 300 &     &      &    \\ \hline
        &       & 311 & 021 & 100 & 230 &      &   \\ \hline
        &       &     & 301 & 031  & 110 & 220   & \\ \hline
        &       &     &    & 231 & 101  & 120  & 210 \\ \hline
    200    &       &     &    &    & 221 & 111  &  130\\ \hline
    \end{array}.$$
    Note that every element of $\Z_4\oplus\Z_4\oplus\Z_2$ appears exactly  once in the array. Moreover, each row and each column contains exactly four filled cells and the sum of the elements in each row and in each column is equal to $([2]_4,[2]_4,[0]_2)$.
\end{ex}

If we do not require that our $\MRS_{\Gamma}(n;2;c)$ is diagonal, we can prove the following. 

\begin{prop}\label{pro:interi pari}
Let $m,n,s,k$ be even integers such that
 $2\leqslant s\leqslant n$, $2\leqslant k\leqslant m$ and $ms=nk$.
Let $c\geqslant 1$. For every abelian group  $\Gamma$ of order $nkc$, there exists an $\MRS_\Gamma(m,n;s,k;c)$.  
\end{prop}

\begin{proof}
Since $\frac{m}{2} \cdot \frac{s}{2}=\frac{n}{2}\cdot \frac{k}{2}$,
by Remark \ref{buchi} we can select, in an array of size $\frac{m}{2}\times \frac{n}{2}$, a set $$\mathcal{C}=\left\{ \left(i_\ell,j_\ell\right): 1 \leqslant \ell \leqslant \frac{nk}{4} \right\}$$ of $\frac{nk}{4}$ cells in such a way that each row and each column contains, respectively, $\frac{s}{2}$ and $\frac{k}{2}$ of these cells.
By Theorem \ref{MRSdia} there exists an  $\MRS_\Gamma\left(2,2;\frac{nkc}{4}\right)$. 
Let $\mathcal{R}=\left\{R_1,\ldots,R_{\frac{nkc}{4}}\right\}$ this $\Gamma$-magic rectangle set, where
$$R_i=\begin{array}{|c|c|} \hline
x_i & y_i \\ \hline
z_i & w_i \\ \hline
\end{array}.$$
Let $\omega$ and $\delta$ be, respectively, the sum of the elements in each row and each column of~$R_i$.

For every $h=1,\ldots,c$, take an empty array $A_h$ of size $m\times n$. 
For every $\ell=1,\ldots,\frac{nk}{4}$,
we fill the four cells
$$(2i_\ell-1, 2j_\ell-1), \;(2i_\ell-1, 2j_\ell),\;(2i_\ell, 2j_\ell-1),\; (2i_\ell, 2j_\ell),$$
with, respectively, the elements  $x_{\ell'}$, $y_{\ell'}$,
$z_{\ell'}$, $w_{\ell'}$ of $R_{\ell'}$,  
where $\ell'=\ell+(h-1)\frac{nk}{4}$.

    Then, each row of  $A_h$ has $2\frac{s}{2}=s$ filled cells and each column has $2\frac{k}{2}=k$ filled cells. So, $A_h$ is a partially filled array where the sum of the entries in each row is equal to  $\frac{s}{2}\omega$ and the sum of the entries in each column is equal to $\frac{k}{2}\delta$.
    Furthermore, by construction of the arrays $R_i$, it is clear that $\displaystyle\bigcup_{i=1}^{c}\E(A_i)=\E(\mathcal{R})$ is equal to $\Gamma$. This proves that the set $\{A_1,\ldots,A_c\}$ is an $\MRS_\Gamma(m,n;s,k;c)$.
\end{proof}

\begin{lem}\label{lem: esistenza con MCD}
Let $m,n,s,k$ be four integers such that
$2\leqslant s\leqslant n$, $2\leqslant k\leqslant m$ and $ms=nk$. Let $c\geqslant 1$ and let $\Gamma$ be an abelian group of order $nkc$.  
If there exists an  $\MRS_\Gamma(k/d, s; md c /k)$, where
$d=\gcd(s,k)$,  then there exists an 
$\MRS_\Gamma(m,n;s,k;c)$.
\end{lem}

\begin{proof}
Let $s_1,k_1,d,e$ be as in \eqref{s1}.
By hypothesis there exists an 
$\MRS_\Gamma(k_1, s_1 d ; e c )$.
Let $R_1,\ldots, R_{e c}$ be the arrays of this set and write $R_a=(r_{i,j}^a)$. 
Denote by $\omega$ and $\delta$, respectively, the sum of the elements of each row and column of $R_a$.
We fill the cells of an $e k_1 \times e s_1$ array $B_b$, where $b=1,\ldots,c$,  as follows:
for every $1\leqslant f \leqslant e$,
every $1\leqslant i \leqslant k_1$ and every
$1\leqslant j \leqslant ds_1 $, 
  we fill the cell 
  $(k_1(f-1)+ i, s_1(f-1)+j)$ of $B_b$
with the element $r_{i,j}^{(b-1)e+f }$, working modulo
$n=e s_1$ with residues in $\{1,\ldots,n\}$ for column indices.

Hence, we obtain a set $\{B_1,\ldots,B_c\}$ of arrays where the sum of the entries of each row is equal to $\omega$, while the sum of the entries of each column is equal to $d\delta$. Furthermore, by construction of the arrays $B_b$, it is clear that 
$\displaystyle\bigcup_{b=1}^{c}\E(B_b)=\Gamma$. This proves that the set $\{B_1,\ldots,B_c\}$ is an $\MRS_\Gamma(m,n;s,k;c)$.
\end{proof}

\begin{ex}
We construct an $\MRS_{\Z_8\oplus\Z_3}(6,3;2,4;2)$. Following the proof of the previous lemma, 
    set $d=\gcd(4,2)=2$, $s_1=1$,  $k_1=2$
    and $e=3$.
    Therefore, taking an  
    $\MRS_{\Z_8\oplus\Z_3}(2,2;6)$  we can construct the set $\MRS_{\Z_8\oplus\Z_3}(6,3;2,4;2)$, whose elements $B_1$ and $B_2$ are given below.
    $$B_1=\begin{array}{|c|c|c|}\hline
    ([0]_8,[0]_3) & ([7]_8,[0]_3) & \\ \hline
    ([3]_8,[0]_3) & ([4]_8,[0]_3) & \\ \hline
     & ([1]_8,[0]_3) & ([6]_8,[0]_3) \\ \hline
     & ([2]_8,[0]_3) & ([5]_8,[0]_3) \\ \hline
    ([7]_8,[2]_3) & & ([0]_8,[1]_3) \\ \hline
    ([4]_8,[1]_3) & & ([3]_8,[2]_3) \\ \hline
    \end{array}, \; B_2=\begin{array}{|c|c|c|}\hline
     ([1]_8,[1]_3) & ([6]_8,[2]_3) & \\ \hline
     ([2]_8,[2]_3) & ([5]_8,[1]_3) & \\ \hline
      & ([0]_8,[2]_3) & ([7]_8,[1]_3) \\ \hline
      & ([3]_8,[1]_3) & ([4]_8,[2]_3) \\ \hline
    ([6]_8,[1]_3) & & ([1]_8,[2]_3) \\ \hline
    ([5]_8,[2]_3) & & ([2]_8,[1]_3) \\ \hline
    \end{array}.$$
    Note that every element of $\Z_8\oplus\Z_3$ appears  exactly once and in a unique array. Moreover, each row has $2$ filled cells,  each column has $4$ filled cells, the sum of the entries in each row is equal to $([7]_8,[0]_3)$ and the sum of the entries in each column is equal to 
    $([6]_8,[0]_3)$.
\end{ex}

\begin{cor}\label{cor123}
Let $m,n,s,k$ be four integers such that
$2\leqslant s\leqslant n$, $2\leqslant k\leqslant m$ and $ms=nk$.
Let $c\geqslant 1$ and let $\Gamma\in \Upsilon$ be an abelian group of order  $nkc$.
There exists an $\MRS_\Gamma(m,n;s,k;c)$ in each of the following three cases:
\begin{itemize}
\item[(1)] $s$ and $k$ are two distinct odd integers;
\item[(2)] $s$ and $k$ are two distinct  integers such that $s\equiv k\equiv 2 \pmod 4$ and $s,k>2$;
\item[(3)] $\frac{k}{\gcd(s,k)}$ is an even integer which is not a $2$-power, while $s$ is an odd integer.
\end{itemize}
\end{cor}

\begin{proof}
Let $s_1,k_1,d,e$ be as in \eqref{s1}.
In cases (1) and (2), since $s\neq k$, up to transposition,  we can assume $k_1\geqslant 3$. 
In case (3), since $k_1$ is not a
$2$-power, we have $k_1\geqslant 6$.
In all three cases, since $\Gamma \in \Upsilon$, there exists an $\MRS_\Gamma(k_1,s; ce  )$
by Theorem \ref{th:esistenza}. 
Hence, the existence of an $\MRS_\Gamma(m,n;s,k;c)$ follows from Lemma 
\ref{lem: esistenza con MCD}.
\end{proof}

\begin{prop}
Let $m,n,s,k$ be four integers such that
$2\leqslant s\leqslant n$, $2\leqslant k\leqslant m$ and $ms=nk$.
Let $c\geqslant 1$ and let $\Gamma$ be an abelian group of order $nkc$.
Suppose that $s$ and $k$ are two distinct odd integers.
There exists an $\MRS_\Gamma(m,n;s,k;c)$
if and only if $\Gamma\in \Upsilon$.
\end{prop}

\begin{proof}
If $\Gamma\in \Upsilon$, then the result follows from Corollary  \ref{cor123}.
If $\Gamma \not \in \Upsilon$, 
then $\Gamma$ has a unique involution and the result 
follows from Corollary~\ref{skodd1}.
\end{proof}

Our main aim would be to provide necessary and sufficient conditions for the existence of an $\MRS_\Gamma(m,n;s,k;c)$: at the moment, we are quite far from achieving this goal.
We conclude this paper with the following result, that almost solves this problem when the two parameters $s$ and $k$ are both even.

\begin{prop}
Let $m,n,s,k$ be four integers such that
$2\leqslant s\leqslant n$, $2\leqslant k\leqslant m$ and $ms=nk$.
Let $c\geqslant 1$ and let $\Gamma$ be an abelian group of order  $nkc$. There exists an $\MRS_\Gamma(m,n;s,k;c)$ in each of the following cases:
\begin{itemize}
\item[$(1)$] $s\equiv k\equiv 0 \pmod 4$;
\item[$(2)$] $s\equiv 2\pmod 4$ and $k\equiv 0\pmod 4$;
\item[$(3)$] $s\equiv 0\pmod 4$ and $k\equiv 2\pmod 4$;
\item[$(4)$] $s\equiv k\equiv 2 \pmod 4$ and $m\equiv n\equiv 0 \pmod 2$.
\end{itemize}
\end{prop}

\begin{proof}
Let $s_1,k_1,d,e$ be as in \eqref{s1}.

\noindent (1) In this case, we have $d\equiv 0\pmod 4$.
By Proposition \ref{diag4} we can construct a diagonal $\MRS_\Gamma(nk_1; d ; c)$:  the existence of an $\MRS_\Gamma(m,n;s,k;c)$ follows from Theorem \ref{sq->rt}.

\noindent (2) In this case, we have $d\equiv 2 \pmod 4$,
$s_1$ is odd while $k_1$ is even.
By Proposition \ref{pro:interi pari} 
we can construct an $\MRS_\Gamma(k_1, s ;  ce)$: the existence of an $\MRS_\Gamma(m,n;s,k;c)$ follows from 
Lemma~\ref{lem: esistenza con MCD}.

\noindent (3) This case follows from (2)  and Remark \ref{tr}.

\noindent (4) By Proposition \ref{pro:interi pari} 
we can construct an $\MRS_\Gamma(m,n; s,k; c)$.
\end{proof}

We remark that, when  $s\equiv k \equiv 2\pmod 4$, the parameter $m$ is even if and only if  $n$ is even. So, to solve the existence problem of an $\MRS_\Gamma(m,n;s,k;c)$ when $s$ and $k$ are both even, what is left open is the case when $s\equiv k \equiv 2 \pmod 4$ and $m\equiv n \equiv 1 \pmod 2$.

\section*{Acknowledgments}

The results of this paper are part of  third author's work for her master's thesis.
The first and the second authors are partially supported by INdAM-GNSAGA.

\end{document}